\documentclass[final, onefignum,onetabnum]{siamart190516}

\usepackage{lipsum}
\usepackage{amsfonts}
\usepackage{graphicx}
\usepackage{epstopdf}
\ifpdf
  \DeclareGraphicsExtensions{.eps,.pdf,.png,.jpg}
\else
  \DeclareGraphicsExtensions{.eps}
\fi

\newsiamremark{remark}{Remark}
\newsiamremark{hypothesis}{Hypothesis}
\crefname{hypothesis}{Hypothesis}{Hypotheses}
\newsiamthm{claim}{Claim}

\headers{}{}

\title{DeepOnet Based Preconditioning Strategies For Solving Parametric Linear Systems of Equations}

\author{Alena Kopani\v{c}\'akov\'a \thanks{Division of Applied Mathematics, Brown University, USA
   (\email{alena.kopanicakova@usi.ch}).}
\and George Em Karniadakis \thanks{Division of Applied Mathematics, Brown University, USA
   (\email{george\_karniadakis@brown.edu}).} }

\usepackage{amsopn}

 \usepackage[T1]{fontenc}
\usepackage[latin1]{inputenc}
\usepackage{caption}
\usepackage{algorithm,algpseudocode}
\usepackage{graphicx,url}
\usepackage{amsmath,amssymb,amsopn}
\usepackage{tabularx}
\usepackage{color}
\usepackage{calc}
\usepackage{caption}
\usepackage{longtable}
\usepackage{multirow}
\usepackage{mathrsfs}
\usepackage{array}
\usepackage{hyperref}
\usepackage{tikz}
\usepackage{booktabs}
\usepackage{pgfplotstable}
\usepackage{pgfplots}
\usepackage[outline]{contour}
\usepackage{calrsfs}
\usepackage{graphicx}
\usepgfplotslibrary{groupplots}
\usetikzlibrary{matrix}
\usetikzlibrary{automata, external, shapes.geometric, fit}
\pgfplotsset{compat=1.13}

\def \ev{\vec{e}}
\def \fv{\vec{f}}

\def \pv{\vec{p}}
\def \qv{\vec{q}}
\def \rv{\vec{r}}

\def \uv{\vec{u}}
\def \vv{\vec{v}}
\def \wv{\vec{w}}
\def \xv{\vec{x}}
\def \yv{\vec{y}}
\def \zv{\vec{z}}

\def \mm{\mat{m}}

\def \rm{\mat{r}}

\def \ym{\mat{y}}

\def \Am{\mat{A}}

\def \Cm{\mat{A_c}}
\def \Dm{\mat{D}}
\def \Em{\mat{E}}

\def \Hm{\mat{H}}
\def \Im{\mat{I}}

\def \Mm{\mat{M}}

\def \Qm{\mat{C}}
\def \QQm{\mat{Q}}
\def \Pm{\mat{P}}
\def \Rm{\mat{R}}

\def \Vm{\mat{V}}

\def \Zm{\mat{Z}}
\def \Id{\mat{I}}

\def \R{\mathbb{R}}			\def \N {\mathbb{N}}			\def \E {\mathbb{E}}

\DeclareMathAlphabet{\pazocal}{OMS}{zplm}{m}{n}

\renewcommand{\vec}[1]{\boldsymbol{#1}}
\newcommand{\mat}[1]{\boldsymbol{{#1}}}

\algnewcommand\algorithmiconput{\textbf{Constants:}}
\algnewcommand\algorithmicinput{\textbf{Input:}}
\algnewcommand\algorithmicoutput{\textbf{Output:}}
\algnewcommand{\algorithmicgoto}{\textbf{go to}}

\algnewcommand\Constants{\item[\algorithmiconput]}
\algnewcommand\Input{\item[\algorithmicinput]}\algnewcommand\Output{\item[\algorithmicoutput]}\algnewcommand{\Goto}[1]{\algorithmicgoto~\ref{#1}}

\definecolor{myblack}{RGB}{53, 53, 53}
\definecolor{myblue}{RGB}{40, 75, 99}
\definecolor{myred}{RGB}{192, 50, 33}
\definecolor{myyellow}{RGB}{255, 166, 48}
\definecolor{mywhite}{RGB}{240, 237, 238}
\definecolor{mygreen}{RGB}{0, 102, 0}

\definecolor{green1}{RGB}{9, 82, 86}
\definecolor{green2}{RGB}{8, 127, 140}
\definecolor{green3}{RGB}{6, 167, 125}
\definecolor{green4}{RGB}{79, 109, 122}
\definecolor{green5}{RGB}{192, 214, 223}
\definecolor{violet}{RGB}{26,69,131}

\definecolor{checkgreen}{rgb}{0,0.6,0}
\definecolor{phase1}{rgb}{0.008,0.655,1.000}
\definecolor{phase2}{rgb}{0.016,0.75,0.700}
\definecolor{phase3}{rgb}{0.929,0.35,0.700}
\definecolor{icsyellow}{cmyk}{0.00,0.11,0.53,0.00}

\definecolor{blackmy}{RGB}{38, 70, 83}
\definecolor{bluemy}{RGB}{39, 125, 161}
\definecolor{greenmy}{RGB}{42, 167, 143}
\definecolor{yellowmy}{RGB}{233, 196, 106}
\definecolor{brownmy}{RGB}{244, 162, 97}
\definecolor{redmy}{RGB}{249, 65, 68}

\definecolor{darkbluemy}{RGB}{65, 59, 147}
\definecolor{lightbluemy}{RGB}{71, 139, 194}
\definecolor{greenmy}{RGB}{98, 173, 153}
\definecolor{darkorangemy}{RGB}{230, 142, 52}
\definecolor{lightorangemy}{RGB}{217, 172, 59}

\definecolor{blue1}{RGB}{1, 58, 99}
\definecolor{blue2}{RGB}{42, 111, 151}
\definecolor{blue3}{RGB}{70, 143, 175}
\definecolor{blue4}{RGB}{137, 194, 217}

\definecolor{red1}{RGB}{204, 68, 75}
\definecolor{red2}{RGB}{218, 85, 82}
\definecolor{red3}{RGB}{227, 150, 149}
\definecolor{red4}{RGB}{228, 190, 171}

\definecolor{brown1}{RGB}{92,178,112}
\definecolor{brown2}{RGB}{130,194,110}		
\definecolor{brown3}{RGB}{163, 193, 173}

 \usepackage[normalem]{ulem}
\usepackage{soul,xcolor}
\usepackage[many]{tcolorbox}
\usepackage{mathtools}
\usepackage{algorithm}
\usepackage{enumitem}
\usepackage{listings}
\usepackage{relsize}
\usepackage{nameref}
\newcommand{\widebar}[1]{\mkern 1.5mu\overline{\mkern-1.5mu#1\mkern-1.5mu}\mkern 1.5mu}

\usepackage{longtable}

\newcommand{\algorithmiccommentMine}[1]{\bgroup\hfill$\triangleright$~{#1}\egroup}
\newcommand\COMMENTmine[1]{\algorithmiccommentMine{#1}}

\newcommand\oldtext[1]{}
\newcommand\cancel[1]{}

\newcommand\oldtextt[1]{}

\newlist{myitemize}{itemize}{3}
\setlist[myitemize,1]{label=\textbullet,leftmargin=1in}
\setlist[myitemize,2]{label=$\rightarrow$,leftmargin=1em}
\setlist[myitemize,3]{label=$\diamond$}

\usetikzlibrary{spy}
\usepackage{pgfplots}
\usetikzlibrary{intersections}
\usepgfplotslibrary{fillbetween}

\ifpdf
\hypersetup{
  pdftitle={DeepOnet Preconditioners},
  pdfauthor={}
}
\fi

\definecolor{darkbluemy}{RGB}{65, 59, 147}
\definecolor{lightbluemy}{RGB}{71, 139, 194}
\definecolor{greenmy}{RGB}{98, 173, 153}
\definecolor{darkorangemy}{RGB}{230, 142, 52}
\definecolor{lightorangemy}{RGB}{217, 172, 59}

\definecolor{blue1}{RGB}{1, 58, 99}
\definecolor{blue2}{RGB}{42, 111, 151}
\definecolor{blue3}{RGB}{70, 143, 175}
\definecolor{blue4}{RGB}{137, 194, 217}

\definecolor{red1}{RGB}{204, 68, 75}
\definecolor{red2}{RGB}{218, 85, 82}
\definecolor{red3}{RGB}{227, 150, 149}
\definecolor{red4}{RGB}{228, 190, 171}

\definecolor{brown1}{RGB}{92,178,112}
\definecolor{brown2}{RGB}{130,194,110}		
\definecolor{brown3}{RGB}{163, 193, 173}

\usepackage{array}
\newcommand{\PreserveBackslash}[1]{\let\temp=\\#1\let\\=\temp}
\newcolumntype{C}[1]{>{\PreserveBackslash\centering}p{#1}}
\newcolumntype{R}[1]{>{\PreserveBackslash\raggedleft}p{#1}}
\newcolumntype{L}[1]{>{\PreserveBackslash\raggedright}p{#1}}

\DeclareMathAlphabet\bpazocal{OMS}{cmsy}{b}{n}
\usetikzlibrary{patterns}
\begin{document}

\setlength{\belowcaptionskip}{-10pt}
\captionsetup{belowskip=-10pt}
\setlength{\abovedisplayskip}{3.5pt}
\setlength{\belowdisplayskip}{3.5pt}

\maketitle

\begin{abstract}
We introduce a new class of hybrid preconditioners for solving parametric linear systems of equations. 
The proposed preconditioners are constructed by hybridizing the deep operator network, namely DeepONet, with standard iterative methods. 
Exploiting the spectral bias, DeepONet-based components are harnessed to address low-frequency error components, while conventional iterative methods are employed to mitigate high-frequency error components.
Our preconditioning framework comprises two distinct hybridization approaches: direct preconditioning (DP) and trunk basis (TB) approaches. 
In the DP approach, DeepONet is used to approximate an action of an inverse operator to a vector during each preconditioning step. 
In contrast, the TB approach extracts basis functions from the trained DeepONet to construct a map to a smaller subspace, in which the low-frequency component of the error can be effectively eliminated. 
Our numerical results demonstrate that utilizing the TB approach enhances the convergence of Krylov methods by a large margin compared to standard non-hybrid preconditioning strategies. 
Moreover, the proposed hybrid preconditioners exhibit robustness across a wide range of model parameters and problem resolutions.
\end{abstract}



\begin{sloppypar}
\section{Introduction}
Partial differential equations (PDEs) are useful in describing physical phenomena and have become ubiquitous in various scientific and engineering fields.
The solution of PDEs is typically obtained numerically using a discretization technique, e.g., finite difference or finite element (FE) methods.
The discretization process produces a discrete system of equations, which must be solved to obtain the PDE's solution.
Solving such systems of equations in a robust and efficient manner is a long-standing challenge in scientific computing. 
The need for an effective solution strategy becomes even more prevalent in multi-query applications, like uncertainty quantification or in control problems, where the underlying PDE must be solved repeatedly for different parameters, e.g., different material properties, force terms, or boundary conditions.
In these particular scenarios, it is paramount to design a robust and scalable strategy for a wide range of parameters and conditions.

Various approaches for solving large-scale linear systems have been proposed in the literature.
For instance, sparse direct solvers~\cite{pardiso204, amestoy2000mumps}  are often employed for small-scale problems due to their robustness and computationally efficient software implementation. 
However, the computational complexity of direct methods prohibits their applicability to large-scale problems as they do not scale well. 
For example, the LU factorization of sparse linear systems requires $\pazocal{O}(n^{3/2})$ flops in two spatial dimensions and
$\pazocal{O}(n^2)$ flops in three spatial dimensions~\cite{bilgen2017phase}. 
As an alternative, Krylov iterative methods~\cite{saad2003iterative}, such as a Conjugate Gradient (CG) or a Generalized Minimal Residual (GMRES), are commonly employed in practice. 
These methods have a favorable computational cost per iteration and are well-suited for massively parallel computing environments, as their primary building blocks, e.g., matrix-vector or dot product, can be efficiently implemented in parallel. 

However, the convergence speed of Krylov methods deteriorates with the increasing condition number of the linear system.
In the context of discretized PDEs considered in this work, the condition number is known to grow with decreasing mesh size $h$.
For instance, in the case of fourth-order PDEs, the condition number increases as $\pazocal{O}(h^{-4})$.
The possible remedy to improve the condition number of the system and, in turn, to enhance the convergence of the Krylov methods is to employ a preconditioning strategy~\cite{saad2003iterative}.
Popular preconditioners include for example the stationary methods, such as Jacobi, Gauss-Seidel~\cite{saad2003iterative}, but also incomplete factorizations~\cite{manteuffel1980incomplete}, sparse approximate inverse (SPAI) methods~\cite{benzi1999comparative}, deflation techniques~\cite{erhel1996restarted}, multigrid (MG)~\cite{trottenberg2000multigrid} and domain-decomposition (DD) methods~\cite{toselli2004domain}.

The MG and DD methods, collectively interpreted as subspace correction methods~\cite{xu1992iterative, tang2009comparison}, are of particular interest for large-scale problems due to their inherent scalability and favorable, often optimal, computational complexity. 
The idea behind these methods is to solve the underlying PDE in a smaller subspace and then use the obtained subspace solution to improve the solution approximation in the full space. 
The subspace correction methods are constructed using two key algorithmic components: subspace mappings, called transfer operators, and subspace solvers. 
A particular choice of these two components greatly influences the applicability and the convergence properties of the resulting iterative methods. 
Indeed, different problem classes often necessitate the use of distinct algorithmic components, the design of which is often challenging and not well-understood in practice. 

In this work, we propose to greatly enhance the convergence of the Krylov methods by devising preconditioning strategies that utilize recent operator learning approaches, namely DeepONet~\cite{lu2021learning, goswami2022physics}. 
In particular, we propose constructing a wide range of preconditioners by hybridizing standard stationary iterative methods with DeepONet-based components.  
Employing hybridization techniques allows us to improve the convergence of standard iterative methods while retaining their accuracy, thus ensuring convergence to user-specified tolerance. 

In the context of linear PDEs, the hybridization of iterative numerical solvers with machine learning (ML) approaches has recently attracted much attention in the literature. 
The ML techniques have been used on various stages of the solution process, for instance, to automate the parameter selection~\cite{arisaka2023principled,hospedales2021meta}, or to obtain a suitable initial guess~\cite{huang2020int, luna2021accelerating, ackmann2020machine}. 
Moreover, several approaches have been proposed for directly intertwining the linear algebra of iterative solvers with ML approaches. 
For example, a hybrid method, tailored explicitly for the Navier-Stokes simulations, has been proposed  in~\cite{tompson2017accelerating}.
This method utilized the standard operator splitting approach but incorporated a convolutional network (CNN) for the pressure projection step, thus approximating the inverse of the discrete Poisson equation required to satisfy the incompressibility constraint. 
The authors of~\cite{hsieh2019learning} have proposed to modify updates made by a numerical solver using a deep neural network (DNN). 
The solver-in-the-loop approach~\cite{um2020solver} was proposed to correct the errors not captured by the discretized PDE using the ML model, which is trained by interacting with the numerical solver. 
A different approach, based on meta-learning of the superstructure of numerical algorithms via recursively recurrent NN, has been proposed in~\cite{doncevic2022recursively}. 

A significant focus has also been given to improving the convergence properties of the Krylov methods. 
For instance, the novel Krylov-subspace-recycling CG method, which employs goal-oriented proper orthogonal decomposition (POD), has been proposed in~\cite{carlberg2016krylov}. 
In \cite{kaneda2022deep}, the authors have employed a CNN to improve the quality of search directions generated by the CG algorithm.
A two-step ML-based approach has been proposed in~\cite{nikolopoulos2022ai}, where the DNN was used to obtain a suitable initial guess and the algebraic MG preconditioner, configured with POD-based transfer operators, has been employed to precondition the Krylov method.
In~\cite{gotz2018machine}, the authors have induced a suitable sparsity pattern for the block-Jacobi preconditioner using CNN. 
Along similar lines, the non-zero pattern for the ILU preconditioner has been predicted using a CNN in~\cite{stanaityte2020ilu}. 
Furthermore, the authors of~\cite{ichimura2020fast}  have proposed to approximate the inverse of linear system by training DNN to approximate the Green's functions. 
A different strategy has been pursued in~\cite{ruelmann2018prospects}, where the authors learn SPAI preconditioners. 
In~\cite{li2023learning}, the graph neural network has been trained to obtain an approximate decomposition of the system matrix, which is then used to precondition the CG method. 

In the realm of preconditioning, researchers have dedicated substantial efforts to the advancement of hybridized subspace correction methods with the goal of not only improving convergence but also enhancing robustness and scalability.
Considering DD methods, ML approaches have been used, for example, to determine overlap or subdomain/interface boundary conditions, c.f.~\cite{burrows2013learning,taghibakhshi2023mg, taghibakhshi2022learning}. 
Alternatively, the physics-based ML surrogates have been utilized to replace the discretization and solution process of the subproblems~\cite{li2020deep,li2019d3m}. 
A particular focus has been to enhancing the construction of the coarse spaces. 
For example, an adaptive coarse space obtained by predicting the necessary location of coarse space constraints in finite element tearing and interconnecting dual-primal (FETI-DP)  and generalized Dryja-Smith-Widlund (GDSW) frameworks have been proposed in~\cite{heinlein2019machine, klawonn2022learning}  and~\cite{heinlein2021combining}, respectively. 
In~\cite{chung2021learning}, the authors have proposed the balancing domain decomposition by constraints (BDDC) method,
which utilizes the adaptive coarse space obtained using DNN. 
A different approach has been pursued in~\cite{ciaramella2022spectral}, where PCA and NNs were used to numerically construct the spectral coarse spaces for sub-structured Schwarz methods.

Many researchers have proposed intertwining MG algorithms with ML approaches. 
For instance, a close connection between linear MG and CNNs has been investigated in~\cite{he2019mgnet, he2023interpretive, alt2023connections}. 
The ML approaches have been used to enhance the performance of existing MG methods by learning their parameters, see for example~\cite{margenberg2022neural, margenberg2021deep, antonietti2021accelerating, katrutsa2020black,grebhahn2016performance}.
Moreover, several authors have taken advantage of ML techniques to enhance the design of transfer operators, see for instance~\cite{weymouth2022data, luz2020learning, greenfeld2019learning, taghibakhshi2021optimization, wang2023learning}.
ML techniques have also been used to design novel smoothers in~\cite{huang2022learning, chen2022meta}. 
In addition, several approaches that take advantage of spectral bias in order to design effective coarse space solvers exist in the literature. 
These approaches utilize the fact that DNNs are very efficient in learning the low-frequency components of the error and, therefore, serve as natural coarse space solvers. 
For instance, the authors of~\cite{cui2022fourier} have proposed a Fourier neural solver, which provides the frequency space correction to eliminate the low-frequency components of the error.
Following a similar line of thought, the MG-based preconditioner for the Helmholtz problem, which utilizes U-Net-based coarse space correction, has been proposed in~\cite{azulay2022multigrid,lerer2023multigrid}.  
Finally, we highlight the operator learning preconditioner for Richardson iteration, which has been proposed in~\cite{zhang2022hybrid}. 
This method, termed HINTS, utilizes the DeepONet~\cite{lu2019DeepONet} to reduce low-frequency components of the error, while high-frequency components of the error are addressed by a standard stationary method.
The HINTS approach has been shown to generalize well for problems with varying geometries in~\cite{kahana2022geometry}.

Motivated by the generalization properties of the DeepONet reported in~\cite{lu2019DeepONet, kahana2022geometry}, we propose a novel class of DeepONet-based preconditioners for Krylov methods. 
Using the subspace correction framework, we discuss how to construct multiplicative and additive preconditioners by effectively intertwining standard iterative methods with DeepONet-based components.
Our hybridization strategy involving DeepONet and iterative methods comprises two distinctive approaches. 
The first approach is motivated by the HINTS framework~\cite{lu2019DeepONet, kahana2022geometry} and 
leverages inference through the DeepONet to approximate the application of the subspace solution operator to a vector. 
The second approach presents a key contribution to this work and utilizes the DeepONet to construct the transfer operators.
These transfer operators are constructed using basis functions extracted from the DeepONet and serve as crucial components in formulating subspace problems.

We demonstrate the effectiveness of the proposed preconditioning framework using several benchmark problems, including standard diffusion, diffusion with jumping coefficients, and the indefinite Helmholtz problem.
Our numerical results illustrate that by employing the proposed DeepONet-based preconditioners, we can effectively enhance the convergence, robustness, and applicability of Krylov methods. 
Moreover, we also show that the proposed preconditioners generalize well and maintain efficiency with respect to the problem's parameters and resolution, frequently outperforming standard numerical approaches by a significant margin.

This paper is organized as follows:  
In \cref{sec:model_problem} we introduce linear parametric PDEs and discusses their high-fidelity and low-fidelity numerical approximations. 
In \cref{sec:numerical_solution}, we review  Krylov methods and discuss the abstract preconditioning framework. 
In \cref{sec:onet_prec} we discuss how to utilize DeepONet to construct efficient preconditioning components.
In \cref{sec:impl}, we describe benchmark problems, which we employ for testing and demonstrating the capabilities of the proposed DeepONet preconditioning framework.
Finally, in \cref{sec:num_results} we demonstrate the numerical performance of the proposed preconditioning framework.
Finally, we include a summary and possible future work in \cref{sec:conclusion}. 

 \section{Parametric linear partial differential equations}
\label{sec:model_problem}
In this work, we aim to solve the parametrized steady-state PDEs.
We consider an input parameter vector~${\boldsymbol{\theta} \in \boldsymbol{\Theta}}$, which contains the problem parameters, e.g., material parameters,  boundary conditions, or source term. 
The parameter space~$\boldsymbol{\Theta}$ is considered to be a closed and bounded subset of $\R^P$, i.e., $\boldsymbol{\Theta} \subset \R^P$, $P \geq 1$.
Let $\Omega \subset \R^{d}$, ${d \in \{1,2,3\}}$ be a computational domain and $\pazocal{V}=\pazocal{V}(\Omega)$ be a suitable Hilbert space, with its dual~$\pazocal{V}'$. 
The parametric problem is given in abstract strong form as: 
For a given~$\boldsymbol{\theta} \in \boldsymbol{\Theta}$, find the solution~$u(\boldsymbol{\theta}) \in \pazocal{V}$ such that
\begin{equation}
\begin{aligned}
\pazocal{A}(\boldsymbol{\theta})  u(\boldsymbol{\theta}) &= f(\boldsymbol{\theta}),  \quad \text{in} \  \pazocal{V}', 
\end{aligned}
\label{eq:param_pde}
\end{equation}
where~$\pazocal{A}(\boldsymbol{\theta}): \pazocal{V} \rightarrow \pazocal{V}'$ denotes a differential operator and~$f(\boldsymbol{\theta})$ is a linear continous form on $\pazocal{V}$ that is an element of $\pazocal{V}'$. 

Problems of this type arise in many applications, such as uncertainty quantification, design optimization, or inverse problems.
The computational cost of these applications can be excessive, especially if the problem parameters fall into a particular range where a high-fidelity solution is required. 
In such instances, the PDE must be solved multiple times with high accuracy.
We aim to accelerate the solution of such problems by leveraging of low-fidelity numerical approximations, namely DeepONet.

\subsection{High-fidelity numerical approximation (FE discretization)} 
In this work, we obtain a high-fidelity solution of~\eqref{eq:param_pde} using the variational principle and the finite element (FE) method. 
Thus, we restate the equation~\eqref{eq:param_pde} to its weak formulation as follows: 
Given parameters~$\boldsymbol{\theta} \in \boldsymbol{\Theta}$, find ${u(\boldsymbol{\theta}) \in \pazocal{V}}$, such that
\begin{equation}
\begin{aligned}
a(u(\boldsymbol{\theta}), v; \boldsymbol{\theta}) &= f(v; \boldsymbol{\theta}), \qquad \text{for all} \ v \in \pazocal{V},
\end{aligned}
\label{eq:weak_form}
\end{equation}
where the parametrized bilinear form  $a(\cdot, \cdot; \boldsymbol{\theta}): \pazocal{V} \times \pazocal{V} \rightarrow \R$ encodes the differential operator and it is obtained from $\pazocal{A}(\boldsymbol{\theta})$ as
\begin{align}
a(u, v; \boldsymbol{\theta}) = \prescript{}{\pazocal{V'}}{\langle} \pazocal{A}(\boldsymbol{\theta}) u, v \rangle_{\pazocal{V}}, \quad \text{for all} \ u, v \in \pazocal{V}.
\end{align}
The linear form $f(\cdot; \boldsymbol{\theta}): \pazocal{V} \rightarrow \R$ is given as ${f(v; \boldsymbol{\theta}) = \prescript{}{\pazocal{V'}}{\langle} f(\boldsymbol{\theta}), v \rangle_{\pazocal{V}}}$. 
Throughout this work, we assume that~$f(\cdot; \boldsymbol{\theta})$ is continuous and linear, and that $a(\cdot, \cdot; \boldsymbol{\theta})$ is continuous and inf-sup stable over~$\pazocal{V} \times \pazocal{V}$, for all~$\boldsymbol{\theta} \in \boldsymbol{\Theta}$.
 Given that those assumptions\footnote{Under stronger assumptions, such as coercivity of~$a(\cdot, \cdot; \boldsymbol{\theta})$, the uniqueness of the solution follows directly from the Lax-Millgram theorem~\cite{quarteroni2014reduced}.} are satisfied then the parametric problem~\eqref{eq:weak_form} is well-posed, i.e., it has a unique solution for every~$\boldsymbol{\theta} \in \boldsymbol{\Theta}$, c.f., Ne\v{c}as theorem~\cite{quarteroni2014reduced}.

In order to solve the problem~\eqref{eq:weak_form} numerically, we consider the mesh $\pazocal{T}$, which encapsulates the domain $\Omega$, and the finite-element space~${\pazocal{V}_h \subset \pazocal{V}}$.
The FE space $\pazocal{V}_h$ is spanned by nodal basis functions~$\{ \psi_i\}_{i=1}^n$, which we use to approximate~$u(\boldsymbol{\theta})$ as follows $u(\boldsymbol{\theta}) \approx u_h({\boldsymbol{\theta}}) = \sum_{i=1}^n \psi_i(x) \uv_i^{\boldsymbol{\theta}}$.
Here, the vector $\uv^{\boldsymbol{\theta}} \in \R^n$ contains the nodal coefficients of the solution. 
Inserting~$u_h({\boldsymbol{\theta}})$ into the weak form~\eqref{eq:weak_form}, we obtain the following discrete problem: 
Given~$\boldsymbol{\theta} \in \boldsymbol{\Theta}$, find $u_h({\boldsymbol{\theta}}) \in \pazocal{V}_h$, such that
\begin{align}
a(u_h(\boldsymbol{\theta}) , v_h; \boldsymbol{\theta}) = f(v_h; \boldsymbol{\theta}), \qquad \text{for all} \ v_h \in \pazocal{V}_h.
\label{eq:weak_discrete_form}
\end{align}
We can find the solution of~\eqref{eq:weak_discrete_form} by solving the following system of linear equations:
\begin{align}
\Am^{\boldsymbol{\theta}}  \uv^{\boldsymbol{\theta}}  = \fv^{\boldsymbol{\theta}}
\label{eq:lin_system_of_eq}
\end{align}
for the nodal coefficients~$\uv^{\boldsymbol{\theta}}  \in \R^n$.
The matrix~$\Am^{\boldsymbol{\theta}}  \in \R^{n \times n}$ and vector~$\fv^{\boldsymbol{\theta}} $ depend affinely on the parameters $\boldsymbol{\theta}$ and their elements are given as 
\begin{align*}
\Am_{ij}^{\boldsymbol{\theta}}  = a(\psi_j, \psi_i; \boldsymbol{\theta}), \qquad \ \  \fv_i^{\boldsymbol{\theta}}  = f(\psi_i; \boldsymbol{\theta}), \qquad \  \ 1 \leq i,j \leq n.  
\end{align*}
In many practical applications, the system~\eqref{eq:lin_system_of_eq} must be assembled and solved for many different parameters~$\boldsymbol{\theta}$. 
This requires significant computational resources, especially when $n$ is very large.
Our goal is to lower these computational resources by utilizing novel, hybrid preconditioning strategies. 
For the remainder of this paper, we omit the superscripts related to the specific parameters $\boldsymbol{\theta}$ in order to simplify the presentation of devised methods.

\subsection{Low-fidelity numerical approximation (DeepONet)}
In order to obtain a low-fidelity solution of~\eqref{eq:param_pde}, we utilize an operator learning approach, namely DeepONet~\cite{lu2021learning}. 
The idea behind the DeepONet is to approximate a mapping between infinite-dimensional function spaces. 
In the context of the parametric PDE considered in this work, this can involve various scenarios. 
For example, we can learn a mapping from a parametrized right-hand side or boundary conditions to the solution of the underlying PDE. 
This section briefly introduces DeepONet, while particular extensions required for constructing efficient preconditioners will be discussed in \cref{sec:onet_prec}.

\subsubsection{Single-input DeepONet}
Let $\pazocal{Y}$ be an infinite-dimensional Banach space.
Our goal is to learn a mapping ${\pazocal{G}: \pazocal{Y} \rightarrow \pazocal{V}}$.
The DeepONet approximates the mapping~$\pazocal{G}$ by using two sub-networks, called branch and trunk.
The branch network ${B: \R^{ny} \rightarrow \R^p}$ encodes the input function, while the trunk network ${T: \R^{d} \rightarrow \R^p}$ is used to encode the domain of the solution.
The branch and trunk networks are then combined to approximate $\pazocal{G}$ as follows
\begin{align}
\pazocal{G}(y)(\boldsymbol{\xi}) \approx \sum^p_{k=1} B_k(\yv) \times T_k(\boldsymbol{\xi}), \quad y \in \pazocal{Y}, \ \boldsymbol{\xi} \in \Omega,
\label{eq:single_input_DeepONet}
\end{align}
where~$\times$ denotes a point-wise multiplication.
The sets $\{B_k\}_{k=1}^p$ and~$\{T_k\}_{k=1}^p$ denote $p$ outputs of the branch and trunk networks, representing the coefficients and the basis functions of the solution's approximation, respectively.

The input to the branch network is represented by the finite-dimensional approximation of the infinite-dimensional input function $y \in \pazocal{Y}$.
In this work, we approximate a function $y$ in a finite-dimensional space $\pazocal{Y}_h$ by evaluating $y$ at  $nb$ points $\{\qv_j\}^{nb}_{j=1}$, called sensor locations, giving rise to a finite-dimensional vector~$\yv \in \R^{nb}$.

\subsubsection{Multiple-input DeepONet}
The PDE given by~\eqref{eq:param_pde} is quite often parametrized, such that multiple input functions have to be considered. 
For example, \eqref{eq:param_pde} might undergo parametrization with respect to the right-hand side as well as the boundary conditions. 
In this case, the DeepONet has to be trained to approximate a map from $nf$ input functions, denoted by~$\{ y^{l} \}_{l=1}^{nf}$, to the solution of PDE, i.e., 
\begin{align}
\pazocal{G}:  \pazocal{Y}^1 \times \cdots \times  \pazocal{Y}^{nf} \rightarrow \pazocal{V}. 
\label{eq:DON_multi_dim}
\end{align}
To this aim, we follow~\cite{jin2022mionet} and extend the architecture of DeepONet to have $nf$ branches. 
The $l$-th branch encodes the $l$-th input function~$y^{l}$, while the trunk still encodes the coordinate point at which the solution is evaluated. 
The nonlinear operator~$\pazocal{G}$ given by~\eqref{eq:DON_multi_dim} is then approximated as follows
\begin{align}
\pazocal{G}(y^1, \ldots, y^{nf})(\boldsymbol{\xi}) \approx 
\sum^p_{k=1} B_k^1(\yv^1)  \times \cdots \times B_k^{nf}(\yv^{nf})  \times T_k(\boldsymbol{\xi}),
\label{eq:multiinput_DeepONet}
\end{align}
where $\boldsymbol{\xi} \in \Omega$. 
The symbol $B_k^l$ denotes the output of the $l$-th branch network and  $\yv^l \in \R^{nb_l}$ stand for the discrete representation of the function~$y^{l}$. 
Note that each input function~$y^{l}$ can be discretized using a different set of sensor locations $\{ \qv_j^l \}_{j=1}^{nb_l}$.

\begin{figure}
\begin{minipage}{0.42\linewidth}
\centering
\vspace{1.45cm}
\rotatebox{-90}{
\scalebox{0.7}{
\includegraphics{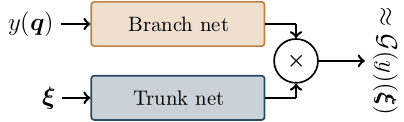}}}
\end{minipage}
\hfill
\begin{minipage}{0.55\linewidth}
\rotatebox{-90}{\scalebox{0.7}{\includegraphics{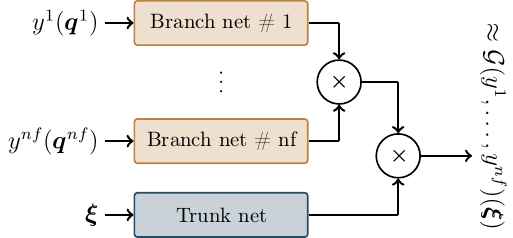}}}
\end{minipage}
\caption{An example of single/multi-input (left/right) DeepONet, see~\cite{lu2022comprehensive}.}
\label{fig:DeepONets}
\end{figure}

An illustration of single and multi-input DeepONet architectures is depicted in \cref{fig:DeepONets}. 
Note that the branch and trunk networks can be represented by different types of DNNs, e.g., feed-forward networks (FFN), or convolutional networks (Conv).

\subsubsection{Training}
The optimal parameters of the DeepONet are found by training. 
For this purpose, we utilize the dataset~$\pazocal{D}$ of $N_S$ samples, given as
\begin{align}
\pazocal{D}=\{(\yv_j^1, \yv_j^2, \ldots, \yv_j^{nf},  \bar{\boldsymbol{\xi}}_j,  \uv_j) \}_{j=1}^{N_s},
\end{align}
which takes into account the discretized functions~$\{ y^{i} \}_{i=1}^{nf}$, denoted as $\{ \yv^{i} \}_{i=1}^{nf}$, where each~${\yv^{i} \in \R^{nb_i}}$. 
Moreover, each sample also contains a set of nodal points, represented by a tensor~${\bar{\boldsymbol{\xi}}_j = [\boldsymbol{\xi}_{j, 1}, \ldots, \boldsymbol{\xi}_{j, n_{\text{don}}}]^{\top} \in \R^{n_{\text{don}} \times d}}$.
The target solution, denoted by $\uv_j \in \R^{n_{\text{don}}}$, represents an approximation of $\pazocal{G}(y^1_j, \ldots, y^{nf}_j)$ at all points given in the set~$\{ \boldsymbol{\xi}_{j, i} \}_{i=1}^{n_{\text{don}}}$.

The training is performed by minimizing the relative error between the output of the DeepOnet and the target solution as
\begin{align}
\mathlarger{\min}_{\wv \in \R^{np}} \ \ 	{\mathlarger{\sum}}_{j=1}^{N_s}  \frac{\Big{\|}\ \bigg( \mathlarger{\sum}^p_{k=1} \big( \prod_{l=1}^{nf} B_l^1(\wv; \yv^l) \big) \times T_k(\wv; \boldsymbol{\xi})  \bigg) - \uv_j \ \Big{\|}^2_2}{\|  \uv_j {\|}^2_2},
\label{eq:training_formula}
\end{align}
where~$\wv$ denotes collectively all parameters of DeepONet.
The symbol~$\prod$ in~\eqref{eq:training_formula} stands for the point-wise product. 
Note that for the clarity of the presentation, the descriptions of single and multi-input DeepONets given in previous sections avoid explicit dependence on the parameters~$\wv$.

 \section{Numerical solution of high-fidelity linear systems using preconditioned Krylov methods}
\label{sec:numerical_solution}
Krylov's methods are considered to be amongst the most efficient iterative solution strategies for solving linear systems of equations. 
Given an initial guess~$\uv^{(0)}$, Krylov methods seek an approximate solution~$\uv^{(i)}$ of~\eqref{eq:lin_system_of_eq} in subspace $\uv^{(0)} + \pazocal{K}_i(\Am, \rv^{(0)})$, where $\pazocal{K}_i(\Am, \rv^{(0)})$ is the Krylov subspace, defined as 
\begin{align}
\pazocal{K}_i(\Am, \rv^{(0)}) := \text{span} \  \{ \rv^{(0)}, \Am \rv^{(0)} \ldots, \Am^{(i-1)} \rv^{(0)} \}. 
\end{align}
Moreover, the residual~$\rv^{(i)}$ is required to be orthogonal to a subspace $\pazocal{L}_i (\Am, \uv^{(0)})$ of dimension $i$. 
Thus, we have to ensure the following Petrov-Galerkin condition~\cite{saad2003iterative} 
\begin{align}
\fv - \Am \uv^{(i)} \perp \pazocal{L}_i (\Am, \uv^{(0)}).
\label{eq:orthogonality_condition}
\end{align}

Existent literature on Krylov methods contains several variants, which differ in a choice of subspace~$\pazocal{L}_i$.
For example, GMRES method~\cite{saad2003iterative} is derived by setting ${\pazocal{L}_i (\Am, \uv^{(0)}) := \Am \pazocal{K}_i(\Am, \rv^{(0)})}$, while CG method is derived by choosing~$\pazocal{L}_i := \pazocal{K}_i$.
The algorithmic details regarding both methods can be found in \ref{sec:krylov_examples}.

\subsection{Preconditioning}
\label{sec:preconditioning}
The convergence rate of Krylov methods depends on two main factors: the condition number of the matrix~$\Am$ and how the eigenvalues of~$\Am$ are clustered.
In practice, preconditioning techniques,  
which transfer the original system into a new one,
can be utilized to improve the efficiency of the Krylov methods. 
To this aim, let us define a non-singular\footnote{In case of CG method, the operator~$\Mm$ is also required to be SPD.} operator $\Mm \in \R^{n \times n}$, called preconditioner. 
We apply a preconditioner from the right side, i.e., 
\begin{align}
\Am \Mm \boldsymbol{\vartheta} = \fv, \qquad \text{where}  \ \ \boldsymbol{\vartheta} = \Mm^{-1} \uv.
\label{eq:prec_system}
\end{align}
In contrast to left-preconditioning, using right-preconditioning allows us to employ flexible variants of Krylov methods, which allow for variable preconditioning.

The convergence speed of the preconditioned Krylov method is determined by the quality of the preconditioner~$\Mm$.
We can investigate the convergence properties of~$\Mm$ by examining its error propagation operator~${\Em \in \R^{n \times n}}$, given as
\begin{align}
\Em = \Id - \Am \Mm,
\end{align}
where~$\Id \in \R^{n \times n}$  denotes an identity matrix. 
Let~$\rho(\Em)$ denote a spectral radius of~$\Em$. 
The preconditioner~$\Mm$ gives rise to convergent iteration if and only if~$\rho(\Em) < 1$. 
Moreover, if $\kappa(\Am \Mm ) \ll \kappa(\Am)$, or if the eigenvalues of $\Am \Mm$ are more clustered than eigenvalues of~$\Am$, then the convergence of the Krylov method is expected to improve.

\subsubsection{Main building blocks for preconditioning}
Ideally, the preconditioner $\Mm$ approximates~$\Am^{-1}$ as closely as possible, but at the same time, it is computationally cheap to apply.
Moreover, the preconditioner shall be scalable and maintain robustness and efficiency across a wide range of parameters.

In this work, we propose a novel class of preconditioners by hybridizing DeepONet with standard iterative methods.
The proposed preconditioners are built upon subspace correction framework~\cite{xu1992iterative, tang2009comparison}. 
To this aim, we assume that there exist two transfer operators, namely a restriction operator $\Rm: \R^n \rightarrow \R^k$ and a prolongation operator $\Pm: \R^k \rightarrow \R^n$, that map data from and to a subspace, respectively.
For the purpose of this work, we assume that $\Pm$ has full rank and that $\Rm:= \Pm^T$. 
Moreover, we define the projection operator~$\boldsymbol{\Pi} \in \R^{n \times n}$, an invertible operator~$\Cm \in \R^{k \times k}$ and the matrix~$\Qm \in \R^{n \times n}$ as 
\begin{align}
\boldsymbol{\Pi} := \Id - \Qm \Am, \qquad \Qm := \Pm \Cm^{-1} \Rm, \qquad \Cm := \Rm \Am \Pm.
\label{eq:all_operators}
\end{align}

Using different transfer operators~$\Pm$ and  $\Rm$ gives rise to different~$\boldsymbol{\Pi}, \Qm$  and $\Cm$,
consequently leading to a different type of preconditioner.
For example, in the context of MG~\cite{trottenberg2000multigrid}, we can identify $\Pm$ and $\Rm$ as the standard prolongation/restriction operators. 
The operator $\Cm$ would then correspond to the coarse-level (Galerkin) operator, while the operator~$\boldsymbol{\Pi}$ can be interpreted as an algebraic form of the coarse-level correction step.  
Similarly, in the context of DD methods~\cite{toselli2004domain}, the matrices~$\Rm$ and~$\Pm$ can be seen as restriction and prolongation operators assembled based on subdomains.

\subsection{Multiplicative preconditioning}
Let us consider a preconditioner~$\Mm$ defined by multiplicatively composing preconditioners~$\Mm_s \in \R^{n \times n}$, where ${s=1, \ldots, S}$, 
such that the error propagation matrix $\Em$ is given as 
\begin{align}
\Em = (\Id - \Am \Mm) =  \prod_{s=1}^{S} (\Id - \gamma_s \Am \Mm_s),
\label{eq:multiplicative}
\end{align}
where the scalar~$\gamma_s \in \R$.

Several well-known preconditioners are of a multiplicative nature, e.g.,  Gauss-Seidel, or multiplicative Schwarz method~\cite{griebel1995abstract}. 
The MG method is perhaps the most prominent, as it is known to be an optimal preconditioner for many problems, e.g., elliptic PDEs.
In the two-level settings, its error propagation matrix is given as
\begin{align}
\Em = (\Id - \Mm_1 \Am)  (\Id - \Mm_2 \Am) = (\Id - \Mm_1 \Am) (\Id - \Qm \Am),
\label{eq:2l_mg}
\end{align}
where $\Mm_1$ represents a smoother, while $\Qm$ is associated with coarse-level step.
The role of smoother $\Mm_1$ is to remove the high-frequency components of the error, while~$\Qm$ is responsible for removing the low-frequency components of the error. 
Note that the quality of the resulting coarse-level correction is determined by the transfer operators~$\Rm$ and  $\Pm$.
In the case of geometric MG methods~\cite{hackbusch2013multi}, the transfer operators are obtained by coarsening/refining the underlying mesh, which can become tedious for highly unstructured meshes. 
In the case of algebraic MG methods~\cite{briggs2000multigrid}, the transfer operators are assembled by exploiting the connectivity of the matrix~$\Am$. 
In this work, we advocate utilizing the coarse space induced by the DeepONet, which does not require knowledge about a hierarchy of meshes nor connectivity of the matrix~$\Am$.

\subsection{Additive preconditioning}
Let us now consider the additive preconditioner~$\Mm$ , constructed by combining $S$ preconditioners as follows
\begin{align}
\Mm =   \sum_{s=1}^S \gamma_s \Mm_s, \qquad \text{with} \qquad \Em = \bigg( \Id - \sum_{s=1}^S \gamma_s \Mm_s \Am \bigg),
\label{eq:additive}
\end{align}
where $\Mm_s \in \R^{n \times n}$ and~$\gamma_s \in \R$, for $s = 1, \ldots, S$. 

Many well-known preconditioners are additive, e.g.,  BDDC, FETI or standard Additive Schwarz method (ASM)~\cite{toselli2004domain}. 
Here, we focus on ASM, which constructs~$\Mm$ by utilizing a decomposition of the computational domain~$\Omega$, such that 
each~$\Mm_s$ approximates the inverse of the local PDE operator related to a specific part of the domain~$\Omega$. 
Notably, the convergence of the ASM tends to deteriorate with an increasing number of subdomains~$S$. 
To achieve the algorithmic scalability, a coarse space\footnote{In the case of DD methods, the coarse space can also be incorporated using multiplicative or deflation approaches.} has to be incorporated. 
In such case, the preconditioner~$\Mm$ is given as
\begin{align}
\Mm = \sum_{s=0}^{S} \gamma_s \Mm_s  = \gamma_{0}  \Qm + \sum_{s=1}^{S} \gamma_s \Mm_s. 
\label{eq:two_level_additive}
\end{align}
The operator~$\Qm$ associated with a coarse space determines the weak and strong scalability of the resulting preconditioner~$\Mm$. 
In the DD literature, several approaches for the construction of~$\Qm$ have been proposed, e.g., Nicolaides coarse-space or spectral coarse spaces; see~\cite {dolean2015introduction} for more details.
While Nicolaides' approach is simple, its applicability is limited. 
On the other hand, spectral approaches are more applicable in practice, but their setup phase (construction of transfer operators) is computationally tedious,  requiring detailed knowledge of the underlying domain decomposition as well as a solution to local/global eigenvalue problems.
This work aims to simplify and automatize the coarse space construction process by utilizing DeepONet, thereby avoiding the need for constructing and solving costly eigenvalue problems.

 \section{DeepONet enhanced preconditioned Krylov methods}
\label{sec:onet_prec}
In this section, we discuss how to enhance the convergence of the Krylov methods by employing the DeepONet-based preconditioners.
In particular, we construct the preconditioners using the subspace correction methodology discussed in the previous section. 
DeepONet-based preconditioning components will serve the purpose of the coarse-level operator, i.e., they are used to reduce the low-frequency components of the error and to capture the global trend.
This choice arises naturally, as the neural networks are known to suffer from spectral bias~\cite{zhang2022hybrid}, i.e., they learn effectively the low-frequency functions. 

\begin{figure}
\centering
\scalebox{0.8}
{\includegraphics{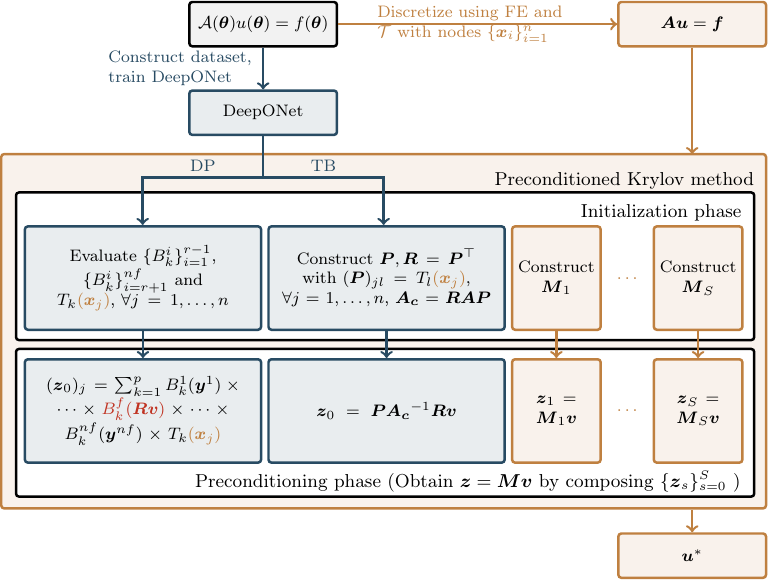}}
\caption{An illustration of a hybrid preconditioning framework.
The DeepONet components are illustrated in blue, while components related to FE discretization and standard iterative methods are illustrated in brown.
The red color depicts the part of DeepONet that must be evaluated during each preconditioning step. }
\label{fig:picture_comparison}
\end{figure}

\subsection{Direct preconditioning (DP) approach}
\label{sec:DP_approach}
One possibility for embedding DeepONet into the subspace correction framework is to approximate the action of~$\Qm$ to a vector directly by means of DeepONet inference. 
In this particular case, the DeepONet has to be designed such that it approximates~$\Am^{-1}\vv$ for a given vector $\vv$. 
To this aim, we have to re-sample the right-hand side of~\eqref{eq:weak_form}, such that the sampling strategy captures distribution of vectors $\vv$ arising during the iterative process.
Moreover, if~\eqref{eq:weak_form} was initially not parametrized with respect to the right-hand side; then the DeepONet architecture has to be adjusted by introducing an additional branch network.

Using the DeepONet, which can approximate the solution of the underlying PDE with respect to different right-hand sides, we can approximate the action of the operator~$\Qm$ to a vector~$\vv$, i.e., $\Qm \vv = \Pm \Cm^{-1} \Rm \vv$, in a component-wise manner as
\begin{align}
(\Qm \vv)_j \approx \label{eq:q_dd_2l_a} 
\sum_{k=1}^p  B_k^1(\yv^1) \times \cdots \times  B_k^{f}(\Rm \vv) \times \cdots \times B_k^{nf}(\yv^{nf}) \times T_k (\boldsymbol{\xv}_j),
\end{align}
for all  $j = 1, \ldots, n$. 
Here, the operator~$\Rm$ maps the vector~$\vv$ to the space~$\pazocal{Y}^f_h$ associated with the sensor locations~$\{\qv^f_i \}_{i=1}^{nb_f}$, where $f$ denotes the index of the branch network used for encoding of the right-hand side.
The practical details regarding the assembly of the operator~$\Rm$ can be found in \ref{sec:appendix_mapping}.

The inference through the branch network in~\eqref{eq:q_dd_2l_a} can be identified with the application of~$\Cm^{-1}$ to $\vv$. 
Notably, the action of the interpolation operator~$\Pm$ is implicitly embedded into the DeepONet architecture by means of inference through the trunk network for all nodes~$\{ \xv_i \}_{i=1}^{n}$, associated with the high-fidelity FE mesh~$\pazocal{T}$.

\subsubsection{HINTS preconditioner}
The HINTS approach proposed in~\cite{zhang2022hybrid} constructs the DeepONet-based preconditioner for a simple Richardson iteration.
In particular, the HINTS preconditioning iteration is given as 
\begin{align}
\uv^{(i+1/2)} &= \uv^{(i)}  + \Mm_1 (\fv - \Am \uv^{(i)}),  \nonumber \\
\uv^{(i+1)}_j = & \uv^{(i+1/2)}_j  +
   \underbrace{\sum_{k=1}^p B_k^1(\yv^1) \times \cdots \times  B_k^{f}\big(\Rm(\rv^{i+1/2})\big) \times \cdots B_k^{nf}(\yv^{nf}) \times T_k (\boldsymbol{\xv}_j)}_{\approx \  \Mm_2 \rv^{i+1/2} \ =  \ \Mm_2 (\fv - \Am \uv^{(i+1/2)})},  \nonumber
\end{align}
for all $j = 1, \ldots, n$. 
The operator~$\Mm_1$ represents a few steps of a standard stationary method, while the action of the operator~$\Mm_2$ is approximated by the DeepONet.

The HINTS iteration  mimics the behavior of the two-level multigrid, as long as the behavior of the operators~$\Mm_1$ and $\Mm_2$ (DeepONet) spectrally complement each other. 
Ideally, $\Mm_1$ eliminates the high-frequency components of the error, while the DeepONet removes the low-frequency components. This can be achieved by constructing dataset~$\pazocal{D}$ such that the right-hand sides are sampled from the space of slowly varying functions.
The authors of HINTS propose  to sample the right-hand sides from the Gaussian random field (GRF) with a relatively large length-scale parameter.

The HINTS method proves effective provided the gap in the spectrum captured by the operator~$\Mm_1$ and DeepONet is not too large. 
However, for large-scale problems, it becomes necessary to introduce additional operators (levels) to address error components in the mid-range frequencies. 
Unfortunately, training DeepONet to eliminate mid-range/high-frequency errors can be quite resource-intensive. 
This is because in order to capture more oscillatory functions, a smaller length-scale parameter has to be used, which, in turn, necessitates discretization of the right-hand side using finer mesh.
Consequently, to take the most advantage of DeepONet, it is more beneficial to identify the components of standard numerical methods responsible for performance loss and focus on designing and training DeepONet to augment that specific aspect.

\subsubsection{From HINTS to preconditioning of F-GMRES}
This work extends the HINTS methodology to the preconditioning of Krylov methods using the DeepONet augmented subspace correction approaches introduced in \cref{sec:preconditioning}. 
It is essential to note that the inference through a DeepONet is a nonlinear operation, which does not preserve properties of the operator~$\Am$, such as linearity, symmetry, or positive definiteness. 
Consequently, the resulting DeepONet preconditioner cannot be used to precondition any Krylov method, but one has to resort to the flexible-GMRES (F-GMRES), which allows for variable preconditioning. 
This might increase computational cost and memory requirements per iteration, but it might also enable overall faster convergence if the DeepONet component proves effective.

\subsection{Trunk basis (TB) approach}
\label{sec:tb}
In this section, we propose a novel approach for constructing DeepONet-based preconditioners.
More precisely, we propose to construct the transfer operators~$\Pm, \Rm$ by extracting the trunk basis functions from the trained DeepONet. 
Each $(j,l)$-th element of the matrix $\Pm \in \R^{n \times k}$, where $k \leq p$, is obtained as follows
\begin{align}
(\Pm)_{jl} = T_l(\boldsymbol{\xv}_j),
\label{eq:matrix_r}
\end{align}
where $T_l(\boldsymbol{\xv}_j) \in \R$ denotes the $l$-th element of the output of the trunk network, evaluated for a given coordinate point~$\boldsymbol{\xv}_j \in \Omega$. 
It is important to note that the trunk basis~(TB) functions can be evaluated at any location~$\xv$, which makes the approach independent of the mesh and discretization strategy used for constructing the dataset~$\pazocal{D}$. 
For instance, we can train DeepONet using a set of $n_{don}$ points $\{ \boldsymbol{\xi}_j \}_{j=1}^{n_{don}} $ associated with a coarse mesh~$\pazocal{T}^c$, while the inference can be performed using $n$ nodes $\{ \boldsymbol{\xv}_j \}_{j=1}^{n} $ associated with a (non-uniformly) refined mesh~$\pazocal{T}$. 

Recalling \cref{sec:preconditioning}, the transfer operators~$\Pm$ and $\Rm := \Pm^{\top}$ are used to construct the coarse space operator~$\Cm \in \R^{k \times k}$ as $\Cm = \Rm \Am \Pm$. 
Thus, the operator~$\Rm$ provides a transformation from a space of dimension $n$ to a space of dimension $k$.
In general, the obtained basis functions have global support, and therefore, operators~$\Rm, \Pm$, and~$\Cm$ are dense. 
The following sections discuss improving the quality, numerical stability, and applicability of the proposed TB approach.

\subsubsection{Incorporating boundary conditions into design of trunk network}
\label{sec:bc}
It is often beneficial to ensure that the coarse level provides corrections only away from the Dirichlet boundary. 
We can ensure that the TB functions do not have support on the Dirichlet boundaries by incorporating the constraints into the architecture of the trunk network. 
Let $\Gamma_D$ be a boundary of~$\Omega$, where the  Dirichlet boundary conditions are imposed. 
Motivated by the truncated basis approach~\cite{kornhuber1997adaptive, kopanivcakova2021multilevel} proposed in the context of MG for bound-constrained optimization, we modify the trunk's output as
\begin{equation}
  \label{eq:mod_bc}
T_k(\xv) \leftarrow   b(\xv) T_k(\xv),
\end{equation}
where the function $b: \R^d \rightarrow \R$ is constructed, such that $b(\xv) = 0$ for all $\xv \in \Gamma_D$ and $b(\xv) \in (0,1]$, otherwise. 
The DeepONet trained using this modified trunk network will, therefore, only learn the solution of the PDE in the interior of the domain, as the products~\eqref{eq:single_input_DeepONet} and~\eqref{eq:multiinput_DeepONet} will be zero at the Dirichlet boundary.

\subsubsection{Orthogonalization of the trunk basis functions}
\label{sec:ortho}
For the coarse operator $\Cm$ to be well-defined, it is crucial to ensure that the transfer operator~$\Pm$ has a full rank.
Moreover, it is often desirable to ensure that the basis functions are orthogonal to each other in order to improve their condition number.
We can ensure that $\Pm$ has full rank and that its columns are orthogonal to each other by performing QR factorization of a tentative~$\widetilde{\Pm}$, obtained using \eqref{eq:matrix_r}, as~$\widetilde{\Pm} = \widebar{\QQm} \widebar{\Rm}$. 
The transfer operator $\Pm$ can be then constructed as a union of columns $\widebar{\QQm}_j$, $j = 1, \ldots,  k$,  of $\widebar{\QQm}$, for which the $j$-th diagonal entry of $\widebar{\Rm}$ satisfies $(\widebar{\Rm})_{jj} \geq \epsilon$, where $\epsilon > 0$.

\begin{figure}
\centering
\scalebox{0.725}{
\includegraphics{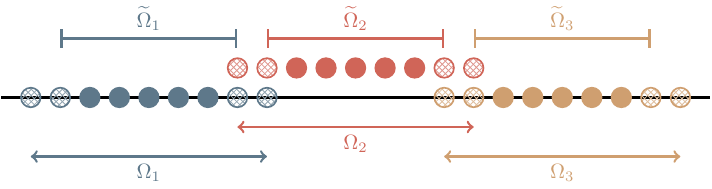}}
\caption{An example of overlapping domain-decomposition of $\Omega$ into three sub-domains $\Omega_1, \Omega_2$ and $\Omega_3$ with overlap = 2.
The hashed pattern depicts the overlapping degrees of freedom.
Symbols $\widetilde{\Omega}_1, \widetilde{\Omega}_2$ and $\widetilde{\Omega}_3$ describe non-overlapping parts of the subdomains.}
\label{fig:dd_picture}
\end{figure}

\subsubsection{Prescribing sparsity pattern}
\label{sec:sparsity_dd}
As the operators~$\Pm$ obtained using~\eqref{eq:matrix_r} is dense, the computational complexity and memory requirements associated with the construction of the coarse space grow significantly with $k$.
To alleviate these requirements, we can prescribe a sparsity pattern during the assembly process.
Following common practice, this can be achieved for example by aggregating subdomain nodes in DD methods or by exploiting the connectivity of the operator~$\Am$ in AMG methods.

In this work, we demonstrate the idea of prescribing the sparsity pattern of $\Pm$  in the context of ASM outlined in~\eqref{eq:two_level_additive}. 
Let the domain $\Omega$ be decomposed into $S$ possibly overlapping sub-domains~$\Omega = \cup_{s=1}^S \Omega_s$. 
Each subdomain is associated with overlapping and non-overlapping index sets, denoted as~$\pazocal{I}_s$ and $\widetilde{\pazocal{I}}_s$, respectively. 
Overlapping index set $\pazocal{I}_s$ contains all degrees of freedom (dofs)  associated with the subdomain~$\Omega_s$. 
To construct non-overlapping index sets, we enforce that $\widetilde{\pazocal{I}}_{i} \cap \widetilde{\pazocal{I}}_{j} = \emptyset$, for $i \neq j$.  
As a consequence, a set of all dofs $\pazocal{I}$ is given as $\pazocal{I} = \cup_{s=1}^S \widetilde{\pazocal{I}}_s = \cup_{s=1}^S {\pazocal{I}}_s$, where~$| \pazocal{I} | := \sum_{s=1}^n n_s = n$, with $n_s := | \widetilde{\pazocal{I}}_s | $. 
An example of such decomposition is depicted in \cref{fig:dd_picture} for an one-dimensional example.

We create an operator $\Pm$ such that it has a block structure using a three-step procedure, motivated by the smoothed aggregation~\cite{vanek1996algebraic}.  
Firstly, we create a dense~$\widetilde{\Pm}$ by extracting $k$ basis functions from the DeepONet, i.e.,
${\widetilde{\Pm} = 
[
\widetilde{\Pm}_1^{\top},  
\widetilde{\Pm}_2^{\top},  
\cdots, 
\widetilde{\Pm}_{S}^{\top} 
]^{\top},
}$
where each block~$\widetilde{\Pm}_{s} \in \R^{n_s \times k }$ is obtained by evaluating trunk basis at a set of points $\{ \boldsymbol{\xv}_j \}_{j \in \widetilde{\pazocal{I}}_s}$, associated with non-overlapping portion of the subdomain $s$. 
Thus $(\widetilde{\Pm}_{s})_{jl}$ is given as $T_l(\boldsymbol{\xv}_j)$, where~$j \in \widetilde{\pazocal{I}}_s$ and ${l \leq k \leq p}$. 
For the purpose of this work, we choose $k$ vectors from $p$ basis functions of trained DeepONet at random.

Secondly, we perform the  QR factorization of each block~$\widetilde{\Pm}_s$ in order to improve the conditioning of the basis functions; c.f., \cref{sec:ortho}. 
The factor~$\widebar{\QQm}_{s}$, where  ${\widetilde{\Pm}_s = \widebar{\QQm}_{s} \widebar{\Rm}_s}$, is then placed to block matrix~$\widebar{\Pm} \in \R^{n \times k \cdot S}$ as follows
\begin{align}
\widebar{\Pm} = 
\begin{bmatrix}
 \widebar{\QQm}_1 & \boldsymbol{0} & \cdots & \boldsymbol{0} \\
\vdots &  \widebar{\QQm}_2& \cdots & \boldsymbol{0} \\
\vdots & \vdots& \ddots & \vdots \\
\boldsymbol{0} & \boldsymbol{0}& \cdots &  \widebar{\QQm}_{S} \\
\end{bmatrix}.
\label{eq:mat_Z}
\end{align}
Note, each $ \widebar{\QQm}_s$ is placed to~$\widebar{\Pm}$ such that its row index set equals to $\widetilde{\pazocal{I}}_s$.

Finally, an optional\footnote{For certain types of problems, such as non-symmetric highly convective problems, it is not advisable to incorporate a prolongation smoothing step.} smoothing step can be performed.
During this step, the smoother is applied to the tentative operator~$\widebar{\Pm}$  in order to ensure that resulting $\pazocal{R}(\Pm^{\top})$ captures the algebraically smooth error more accurately.  
For instance, 
a single iteration of Jacobi prolongation smoothing\footnote{Different transfer operator smoothing approaches might be utilized and give rise to better results.} may be applied to~$\widebar{\Pm}$ as follows
\begin{align}
\Pm  = (\Id - \gamma \Dm^{-1} \Am) \widebar{\Pm},
\end{align}
where~$\gamma \in \R$ and $\Dm = \text{diag}(\Am)$.

The transfer operators~$\Pm, \Rm$ require a storage of~$k$ basis vectors.
However, we will demonstrate numerically in \cref{sec:num_results},  that enforcing the subdomain-based block structure~\eqref{eq:mat_Z} gives rise to an algorithmically scalable ASM even for small values of~$k$. 
Note that the size of~$\Cm \in \R^{(S \cdot k) \times (S \cdot k)}$ grows with the number of subdomains. 
Moreover, the sparsity pattern of the operator $\Cm$ is a result of the sparsity of~$\Rm$.

\subsubsection{Trunk basis versus direct preconditioning approach}
Training the entire DeepONet to utilize only the trunk information might seem wasteful at first. 
However, there are many practical scenarios where training DeepONet is useful for constructing low-fidelity surrogates, such as multi-fidelity uncertainty quantification.  
In such cases, the resulting DeepONet can be used to obtain low-fidelity samples, as well as to expedite the solution process of the high-fidelity problems by providing a suitable initial guess and by constructing efficient preconditioners.

The TB approach differs from the DP approach in several aspects, namely:
\begin{itemize}
\item Using the TB approach, the properties of operator~$\Am$ are preserved. 
For instance, if~$\Am$ is SPD, then also~$\Cm$  will be the SPD matrix. 
Consequently, the resulting DeepONet-based preconditioner can be used to precondition the CG method, while the DP approach would necessitate using the F-GMRES. 
\item The TB approach allows us to use the DeepONet without any modifications to its architecture.
This is especially useful if the problem~\eqref{eq:param_pde} is parametrized such that the right-hand side is kept parameter-free.
A branch network accounting for the right-hand side must be added in such cases.
\item Using the DP approach requires a sampling strategy that effectively captures the distribution of vectors to which the preconditioner shall be applied.
Since the distribution of such vectors is usually unknown, a massive amount of samples is typically required to ensure good performance of the method. 
As a consequence, the construction of the dataset and training of the DeepONet is typically more costly for the DP approach, compared to TB approach.
\end{itemize}

 \section{Benchmark problems and implementation details}
\label{sec:impl}
This section presents a set of benchmark problems, which we employ for testing and demonstrating the capabilities of the proposed DeepONet preconditioning framework. 
Moreover, we provide the details regarding the implementation of our algorithmic framework.

\subsection{Benchmark problems}
We consider benchmark problems formulated in two and three spatial dimensions. 
These problems have been selected either because they are established as standard benchmark problems or due to their documented tendency to pose challenges for the conventional iterative methods.

\begin{figure}
\hspace{0.5cm}
\begin{minipage}{0.3\linewidth}
\scalebox{0.6}{
\includegraphics{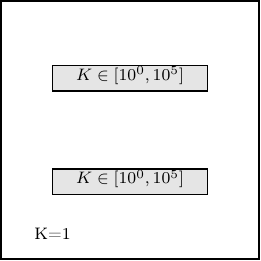}}
\end{minipage}
\begin{minipage}{0.65\linewidth}
\includegraphics[scale=0.216]{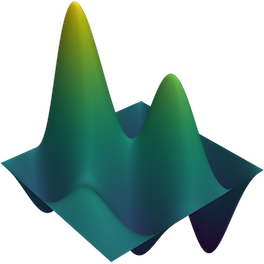}
\includegraphics[scale=0.216]{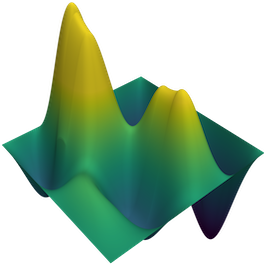}
\includegraphics[scale=0.216]{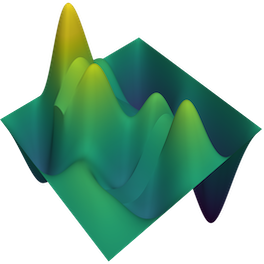}
\includegraphics[scale=0.216]{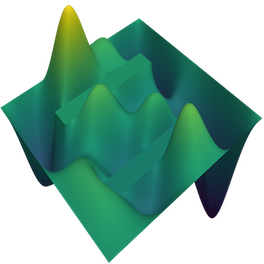}
\end{minipage}
\caption{Left: An illustration of the computational domain with two different channel patterns used for the diffusion equation test with jumping coefficients.
Right: Example of samples used for testing \nameref{sec:jump_diff} example. 
The examples are selected such that the value of~${K \in [1, 10^{5}]}$ increases from left to right.}
\label{fig:jump_solutions}
\end{figure}

\subsubsection[Diff]{Diffusion equation with spatially varying coefficients and forcing term (Diff)}
\label{sec:diff}
We start by considering a diffusion equation given as
\begin{equation}
\begin{aligned}
- \nabla \cdot (K(\xv, \boldsymbol{\theta}) \ \nabla u(\xv) ) &= f(\xv, \boldsymbol{\theta}), \qquad  \  \forall \xv \in \Omega, \\
u(\xv) &= 0, \qquad \qquad \quad  \text{on} \ \partial \Omega,
\label{eq:Laplace}
\end{aligned}
\end{equation}
where~$u$ denotes the solution and $f$ stands for the forcing term.
The equation~\eqref{eq:Laplace} is parametrized in terms of the forcing term and the diffusion coefficient~$K$.

In this first example, we consider~$\Omega := [0,1]^3$, and sample the coefficient~$K$ using Gaussian random fields (GRFs) with mean $\E[K(\xv, \boldsymbol{\theta})]=0.5$ and  the covariances \begin{align}
\text{Cov}(K(\xv,\boldsymbol{\theta}), K( \yv, \boldsymbol{\theta})) = \sigma^2 \exp \bigg(- \frac{\| \xv - \yv \|}{2 \ell^2} \bigg),
\label{eq:sample}
\end{align}
where $\xv$, $\yv$ are two distinct points within the computational domain~$\Omega$. 
The parameters~$\sigma$ and $\ell$ are chosen as $\sigma = 1.0$ and $\ell = 0.1$.
The right hand side~$f$ is also sampled using GRFs, but with $\E[K(\xv, \boldsymbol{\theta})]=0.0$ and the covariance given as in~\eqref{eq:sample}, but with
parameters $\sigma = 1.0$ and $\ell = 0.05$.

\subsubsection[JumpDiff]{Diffusion equation with jumping coefficients (JumpDiff)}
\label{sec:jump_diff}
Next, we consider a scenario in two spatial dimensions, i.e., $\Omega := [0,1]^2$, with fixed right hand side ${f(\xv):=\sin(4 \pi \xv_1) \sin(2 \pi \xv_2) \sin(2 \pi \xv_1 \xv_2)}$ and jumping diffusion coefficients. 
Following~\cite{ciaramella2022spectral}, the diffusion coefficient~$K$ takes on a value one everywhere, except in two rectangular channels, illustrated by grey color in \cref{fig:jump_solutions}. 
In channels, the coefficient $K$ takes on a value from $1$ to $10^5$, which we sample from the distribution~$\log_{10} K \sim \pazocal{U}[0, 5]$.
The FE discretization is performed such that the jumps in the diffusion coefficient are aligned with the edges of elements, see also \cref{fig:jump_solutions}.

\subsubsection[Helm1D]{Helmholtz equation in one spatial dimension (Helm1D)}
\label{sec:helm1D}
Let~${\Omega=[0,1]}$, the Helmholtz equation is used to model the propagation of waves in the frequency domain, and it is given as
\begin{equation}
\begin{aligned}
- \Delta u(\xv) - k_{\text{H}}^2 u(\xv) &= f(\xv, \boldsymbol{\theta}), \qquad   \  \forall \xv \in \Omega, \\
u(\xv) &= 0, \qquad \qquad \ \ \ \text{on} \ \partial \Omega.
\label{eq:helmholz}
\end{aligned}
\end{equation}
We sample the forcing term from GRF with ${\E[K(\xv, \boldsymbol{\theta})]=0.0}$ and the covariance given as in~\eqref{eq:sample}, with~$\sigma = 1.0$ and $\ell = 0.1$. 

The symbol $k_{\text{H}}$ in~\eqref{eq:helmholz} stands for a constant wave number and determines the wavelength~$\lambda$ given as~$\lambda = 2 \pi / k_{\text{H}}$. 
Following common practice in engineering applications~\cite{harari1991finite},  we discretize the problem~\eqref{eq:helmholz} using the FE method, such that the condition $h \leq \frac{\pi}{5 k_{\text{H}}}$ holds. 
We remark that in order to ensure stability, the bound on the term~$h^2 k_{\text{H}}^3$ is also required~\cite{babuska1997pollution}.
Noteworthy, this bound is more restrictive than the condition $h \leq \frac{\pi}{5 k_{\text{H}}}$ for high wave numbers.

\subsubsection[Helm2D]{Helmholtz equation in two spatial dimensions (Helm2D)}
\label{sec:helm2D}
\begin{figure}
\centering
\includegraphics[scale=0.7]{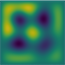}
\hfill
\includegraphics[scale=0.7]{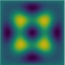} 
\hfill
\includegraphics[scale=0.7]{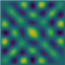}
\hfill
\includegraphics[scale=0.7]{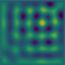} 
\hfill
\includegraphics[scale=0.7]{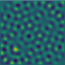}
\hfill
\includegraphics[scale=0.7]{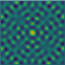}
\caption{
Example of samples used for testing \nameref{sec:helm2D} example. 
The value of~$k_{\text{H}}=15$, $k_{\text{H}}=31$ and $k_{\text{H}}=62$ from left to right, displayed in pairs.}
\label{fig:helm2d}
\end{figure}

We consider Helmholtz equation~\eqref{eq:helmholz} in two spatial dimensions, thus~${\Omega=[0,1]^2}$.
As a forcing term, we consider a Gaussian point source given as~${f(\xv, \boldsymbol{\theta}): = e^{-\frac{1}{2} (\| \xv - \boldsymbol{\theta} \|)/\sigma^2_{\text{H}}}}$, 
where we set~$\sigma_{\text{H}} = 0.8/2^{l-2}$ and $k_{\text{H}} \geq 2^{l-2} \pi / 1.6$.  
Here, the symbol~$l$ denotes the mesh index associated with the FE discretization, as specified in \cref{sec:appB}.
The parameters~${\boldsymbol{\theta} \in \R^2}$ specify the coordinates of the location of the source and are sampled from $\sim \pazocal{U}[0, 1]$. 
Examples of samples used for assessing the performance of proposed preconditioning framework can be found in \cref{fig:helm2d}.

\subsection{Implementation}
In this work, we perform high-fidelity discretization of the benchmark problems using the FE method and the Firedrake framework~\cite{rathgeber2016firedrake}. 
The DeepONets are implemented using PyTorch~\cite{paszke2019pytorch}.
The datasets required for the training of DeepONets are also obtained using the FE method.
Unless specified otherwise, the datasets are constructed by discretization with low-resolution meshes, i.e., $39 \times 39$ elements in 2D and $15 \times 15 \times 15$ elements in 3D.
We train DeepONets in double precision using the  Adam optimizer with a batch size of $1,000$ and learning rate~$10^{-4}$. 
The training terminates as soon as the validation loss does not improve in $10,000$ successive epochs. 
DeepONets are initialized using the Xavier initialization strategy~\cite{kumar2017weight}. 
The details regarding the network architectures, dataset sizes, and training times are summarized in \cref{sec:appB}. 
We point out that our numerical experiments are performed without thorough hyper-parameter tuning, which suggests that the obtained results could be further improved by selecting network architecture and tunning hyper-parameters in a more rigorous way.
Moreover, advanced training methods, such as multilevel~\cite{Kopanicakova_2020c, kopanicakova_22_1, Kopanicakova_2023a} or DD-based~\cite{Kopanicakova_2023b, lee2023training} strategies could be employed to further increase the accuracy and generalization properties of the DeepONets. 

The proposed preconditioning strategies\footnote{The developed code~\cite{precond_git} will be made publicly available upon acceptance of the manuscript.} are implemented using the PETSc library~\cite{balay2019petsc}, which provides an extensive collection of Krylov methods and state-of-the-art preconditioners. 
We implement "\texttt{PC}" components containing DeepONet-based preconditioners via the petsc4py interface. 
This enables the seamless composition of a diverse set of existing preconditioning techniques with a DeepONet component. 
An example of command line options used to compose a DeepONet-based preconditioner can be found in \cref{fig:command_line_options}.

All numerical experiments were conducted on the Piz Daint supercomputer at the Swiss National Supercomputing Centre. 
Each XC50 compute node has an Intel Xeon E5-2690 v3 processor (64 GB) and an NVIDIA Tesla P100 GPU {(16 GB)}.
 \section{Numerical results}
\label{sec:num_results}
In this section, we demonstrate the numerical performance of the proposed DeepONet-based preconditioners. 
We highlight that it is not the purpose of this work to devise the best possible preconditioner for a given problem, but rather, we focus on demonstrating the capabilities and asymptotic convergence properties of the proposed preconditioners. 
The additive and multiplicative preconditioners are examined independently of each other in order to study the impact of the DeepONet-based component on the convergence of the resulting iterative scheme. 
However, all of the proposed techniques could be combined for enhanced convergence.

In order to study the robustness of the preconditioners, their performance is assessed for a wide range of model parameters. 
Moreover, we also study the convergence with respect to the increasing problem size. 
To this aim, the FE discretization is performed using meshes, which are obtained by uniformly refining the coarsest mesh\footnote{The mesh~$\pazocal{T}^1$ is also used to construct the dataset~$\pazocal{D}$ required for training the DeepONet.} $\pazocal{T}^1$ $L-1$ times, giving rise to a hierarchy of $L$ meshes, i.e., $\pazocal{T}^{1}, \ldots, \pazocal{T}^{L}$.  

During all experiments, the preconditioned Krylov methods terminate as soon as the following criteria are satisfied
\begin{align*}
\| \rv^{(i)} \| \leq 10^{-12} \ \ \text{or} \quad \| \uv^{(i)} - \uv^{(i-1)} \|_{\Am} \leq 10^{-12} \ \ \text{or} \quad \frac{\| \rv^{(i)} \|}{\|\rv^{(0)}  \|} \leq 10^{-9}.
\end{align*}
As an initial guess~$\uv^{(0)}$, we always choose a random vector. 
The performance of all considered methods is reported by taking an average over $10$ independent runs, i.e., 
with randomly selected problem parameters and randomly selected initial guesses.

\subsection{Multiplicative preconditioning}
\label{sec:results_multiplicative}
We start our numerical study by investigating the behavior of the multiplicative preconditioners. 

\subsubsection{Convergence of two-level preconditioner for \nameref{sec:diff} problem (DP approach versus TB approach)}
Firstly, we consider a multiplicative preconditioner for the \nameref{sec:diff} example, constructed such that the iteration matrix~$\Em$ is as stated in~\eqref{eq:2l_mg}. 
In order to mimic the behavior of the two-level MG method, we select $\Mm_1$ to be a standard Jacobi with~$\gamma_1 = 2/3$. 
The action of the operator~$\Mm_2:=\Qm$, with  $\gamma_2=1$, is invoked by utilizing DP or TB approaches, as specified in \cref{sec:onet_prec}. 

For the DP approach, we sample the right-hand sides using two distinct strategies.
Firstly, we follow the HINTS methodology~\cite{zhang2022hybrid} and sample right-hand sides from GRF with  ${\E[K(\xv, \boldsymbol{\theta})]=0.0}$ and  the covariance given as in~\eqref{eq:sample}, but with~$\sigma = 1.0$ and $\ell = 0.1$. 
Secondly, we sample the elements of the right-hand side from the Gaussian distribution $\pazocal{N}(0,1)$. 
We denote this strategy as {Rnd}. 
In both cases, the obtained samples are merged with samples used for sampling of the right-hand side, given by the problem definition.
For the TB approach, we consider the operator~$\Rm$ to be a dense matrix constructed by extracting 
$k$ TB functions from the DeepONet. 

{
\begin{table}
\centering
\footnotesize
\tabcolsep=0.145cm
\begin{tabular}{|c|c|l||r|r|r|}
\hline
 \multirow{11}{*}{\rotatebox[origin=c]{0}{$\bpazocal{T}^1$}}  & \multirow{2}{*}{\textbf{KM}}  &\multicolumn{1}{c||}{ \multirow{2}{*}{$\Mm_2= \Qm$ }}   & \multicolumn{3}{c|}{$\boldsymbol{N_S}$} \\ \cline{4-6}
 &  & & $\boldsymbol{2,500}$ & $\boldsymbol{25, 000}$ & $\boldsymbol{100,000}$   \\ \cline{2-6} 
&  \multirow{6}{*}{\rotatebox[origin=c]{90}{F-GMRES(50)}}  & None & \multicolumn{3}{c|}{$\phantom{0}\phantom{0}40.0 \pm 2.6$} \\ \cline{3-6}
 & & DP (GRF) & $35.0 \pm 5.8$ & $32.5 \pm 4.5$ & $29.3 \pm 5.6$ \\  \cline{3-6}
 & & DP ({Rnd}) & $33.1 \pm 4.6$ & $29.4 \pm 4.1$ & $25.1 \pm 3.5$ \\ \cline{3-6}
& & TB (k=8) & $25.4 \pm 1.0$ & $25.3 \pm 0.9$ & $25.5 \pm 1.1$ \\ \cline{3-6}
 & & TB (k=32) & $15.4 \pm 0.4$ & $15.1 \pm 0.3$ & $15.0 \pm 0.3$ \\ \cline{3-6}
 & & TB (k=128) & $8.0 \pm 0.0$ & $8.0 \pm 0.1$ & $8.1 \pm 0.2$ \\ \cline{2-6}  
&  \multirow{4}{*}{\rotatebox[origin=c]{90}{CG}}  & None &  \multicolumn{3}{c|}{$\phantom{0}\phantom{0}36.2 \pm 0.9$}	\\ \cline{3-6}
& & TB (k=8) & $26.1 \pm 0.6$ & $26.0 \pm 0.4$ & $25.9 \pm 0.3$ \\ \cline{3-6}
 & & TB (k=32) & $16.1 \pm 0.3$ & $15.8 \pm 0.1$ & $15.4 \pm 0.1$ \\ \cline{3-6}
 & & TB (k=128) & $8.0 \pm 0.0$ & $8.2 \pm 0.3$ & $8.0 \pm 0.1$ \\ \hline   \hline
  \multirow{13}{*}{\rotatebox[origin=c]{0}{$\bpazocal{T}^4$}} & \multirow{2}{*}{\textbf{KM}}  &\multicolumn{1}{c||}{ \multirow{2}{*}{$\Mm_2 =  \Qm$ }}   & \multicolumn{3}{c|}{$\boldsymbol{N_S}$} \\ \cline{4-6}
 &  & & $\boldsymbol{2,500}$ & $\boldsymbol{25, 000}$ & $\boldsymbol{100,000}$   \\ \cline{2-6}
 &  \multirow{6}{*}{\rotatebox[origin=c]{90}{F-GMRES(50)}}  & None &  \multicolumn{3}{c|}{ $2,252.3 \pm 145.9\phantom{0}$ } \\  \cline{3-6}
 & & DP (GRF) & $564.3 \pm 20.1$ & $549.1 \pm 10.3$ & $500.8 \pm 8.3\phantom{0}$ \\ \cline{3-6}
 & & DP ({Rnd}) & $508.3 \pm 17.3$ & $489.2 \pm 15.7$ & $461.4 \pm 10.6$ \\  \cline{3-6} 
& & TB (k=8) & $557.2 \pm 50.1$ & $552.2 \pm 30.1$ & $522.2 \pm 35.1$ \\ \cline{3-6} 
 & & TB (k=32) & $276.4 \pm 6.4\phantom{0}$ & $262.0 \pm 12.3$ & $258.2 \pm 10.3$ \\ \cline{3-6}
 & & TB (k=128) & $121.8 \pm 2.4\phantom{0}$ & $104.8 \pm 5.2\phantom{0}$ & $104.1 \pm 5.1\phantom{0}$ \\ \cline{2-6}  
&  \multirow{4}{*}{\rotatebox[origin=c]{90}{CG}}    & None & \multicolumn{3}{c|}{$\phantom{0}538.6 \pm 9.3$ } \\ \cline{3-6}
& & TB (k=8) & $393.4 \pm 9.6$ & $380.1 \pm 8.1$ & $375.2 \pm 10.2$ \\ \cline{3-6}
 & & TB (k=32) & $232.0 \pm 9.4$ & $225.1 \pm 5.1$ & $224.2 \pm 5.0\phantom{0}$ \\ \cline{3-6}
 & & TB (k=128) & $123.8 \pm 1.6$ & $117.6 \pm 1.4$ & $117.2 \pm 1.2\phantom{0}$ \\ \hline   
\end{tabular}
\caption{The average number of iterations required by the Krylov method (KM), preconditioned using two level MG, with three pre/post-smoothing steps of Jacobi ($\gamma = 2/3$) and DeepONet-based (DP or TB) coarse step ($\gamma =1$). 
The experiments performed using  \nameref{sec:diff} example, discretized using meshes $\pazocal{T}^{1}$ and $\pazocal{T}^{4}$.}
\label{tab:2l_multiplicative_summary}
\end{table}}

\begin{figure}
\centering
\includegraphics{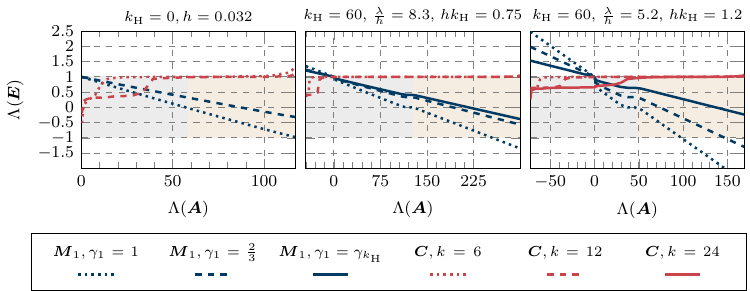}
\caption{
Eigenvalues of the iteration matrix associated with Jacobi smoother $\Mm_1$ with $\gamma_1\in \{1, 2/3, \gamma_{k_{\text{H}}}\}$ and the TB-based operator $\Qm$ with $k \in \{3, 6, 12\}$ plotted against eigenvalues of the associated 1D Helmholtz operator. 
The gray and brown overlays relate to low and high-frequencies parts of the spectrum, respectively.
}
\label{fig:eigen_laplace}
\end{figure}

\Cref{tab:2l_multiplicative_summary} summarizes the obtained results.
As we can see, incorporating a DeepONet-based component helps to reduce the number of required iterations, especially as the problem size increases. 
In addition, we also observe that the obtained results improve when the DeepONet is trained with a larger number of samples.
This is, in particular, true for DP approaches, while the performance of the TB approaches does not vary a lot with an increasing number of samples.
In the context of the DP approach, we also note that the {Rnd} right-hand side sampling strategy yields slightly improved results compared to the GRF sampling strategy. 
These observations confirm our hypothesis that the convergence properties of the composite DP-based preconditioners depend heavily on the choice of sampling strategy and the number of samples used for training the DeepONet.
In the context of the TB approach, we notice that constructing a larger coarse space gives rise to faster convergence.
Notably, faster convergence comes hand in hand with increased memory and iteration costs.

In the end, we also highlight the benefits arising from the fact that TB-based preconditioners can be used to precondition the CG method. 
For instance, for the problem associated with the mesh $\pazocal{T}^4$, the CG with the TB-based preconditioner requires $393.4$ iterations, while $557.2$ iterations are required for the F-GMRES method.
We attribute this difference to the F-GMRES restarts.

We conducted the same series of experiments as described above for all numerical examples outlined in \cref{sec:impl}. 
The results obtained for the TB approach gave rise to conclusions similar to those observed for the \nameref{sec:diff} problem. 
However, using the DP approach did not result in a convergent preconditioner for the \nameref{sec:jump_diff} problem with $K>10^2$ in channels.
We anticipate that different sampling strategies would need to be devised to adequately sample the right-hand sides to which the preconditioner shall be applied. 
Given the superiority and novelty of the TB approach, the following numerical experiments revolve mainly around the TB methodology.

\subsubsection{Numerical analysis of the convergence behaviour of the DeepONet-based composite preconditioner}
\begin{figure}
\centering
\includegraphics{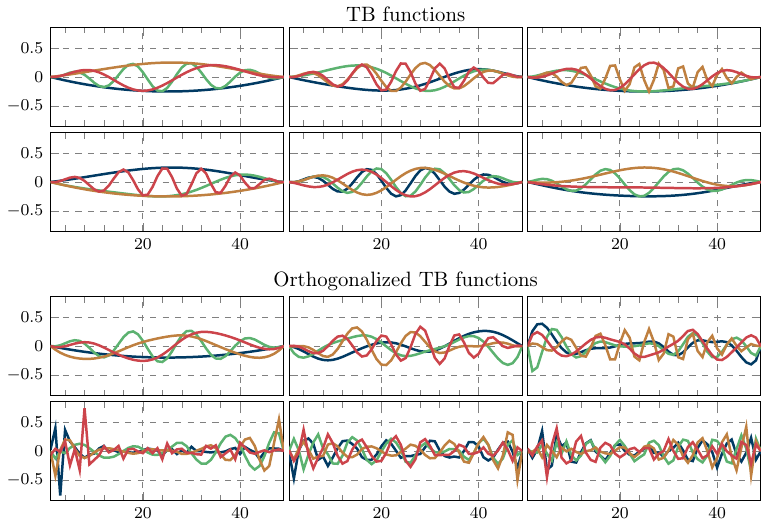}
\caption{Visualization of $24$ randomly selected TB functions extracted from DeepONet ($p=128$), which is trained for \nameref{sec:helm1D} with $k_{\text{H}}=60$ and $h=1/50$, displayed before (top block) and after (bottom block) performing the QR decomposition.}
\label{fig:eigen_fun_helm}
\end{figure}

In this section, we numerically analyze the convergence behavior of the multiplicative preconditioner~\eqref{eq:2l_mg} for the \nameref{sec:helm1D}  example.
To this aim, we study the eigenvalues of the error propagation matrix~$\Em$ associated with $\Mm_1$ and $\Mm_2$.
In place of $\Mm_1$, we again consider Jacobi with different dumping parameters~$\gamma_1$, namely $\gamma_1 \in \{1, 2/3, \gamma_{k_{\text{H}}} \}$. 
Here, the value of $\gamma_{k_{\text{H}}} $ is chosen such that the high-frequency components of the error are reduced the most efficiently for a given $k_{\text{H}}$.
To this aim, we follow~\cite{ernst2011difficult, ernst2013multigrid} and set $\gamma_{k_{\text{H}}} := \frac{2 - k_{\text{H}}^2 h^2}{3 - k_{\text{H}}^2 h^2}$, i.e., $\gamma_{k_{\text{H}}}$ depends on the parameter~$k_{\text{H}}$ and the mesh size~$h$.
The operator, $\Mm_2:=\Qm$, is obtained using the TB approach with an increasing number of TB functions. 
An illustration of the TB functions for this numerical example is shown in \cref{fig:eigen_fun_helm}.
As we can see, performing QR decomposition gives rise to a more diverse set of TB functions, spanning both smooth and highly oscillatory functions. 

\Cref{fig:eigen_laplace} illustrates the obtained results. 
On the left, we depict the results for ${k_{\text{H}}=0}$, i.e., for a Poisson problem. 
As expected, the Jacobi smoother with appropriate damping ($\gamma_1=2/3$) is effective in reducing high-frequency components of the error. 
In contrast, the TB-based coarse space is effective at reducing the low-frequency components of the error. 
Thus, the eigenvalues of error propagation matrix in the lower part of the spectrum (gray region) take on values below one. 
As we can see, the number of reduced low-frequencies depends on the number of extracted TB functions. 
Moreover, we also see that no reduction is achieved for high frequencies, as the associated eigenvalues of the error propagation matrix take on the value one or larger.

Next, we perform the same study with $k_{\text{H}}=60$. 
Here, we utilize two different discretizations, with $81$ (middle part of \cref{fig:eigen_laplace}) and $51$ (right part of \cref{fig:eigen_laplace}) mesh points. 
In the first case, $h=1/80$, the wave-lengths per discretization point ratio is ${\lambda}/{h}=8.3$. 
In the second case, $h=1/50$, the wave-lengths per discretization point ratio is ${\lambda}/{h}=5.3$. 
As we can see from \cref{fig:eigen_laplace}, TB-based coarse space operators behave similarly to what we have observed for the Poisson problem, i.e., the coarse space operates on the lower part of the spectrum. 
Furthermore, we see that the Jacobi method effectively addresses the high-frequency components of the error.
However, we can also notice that at the same time, it amplifies the low-frequencies of the error.
The amplification factor is more prevalent for the case with ${\lambda}/{h}=5.3$ than for the case with ${\lambda}/{h}=8.3$. 
Remarkably, the amplification factor becomes less prevalent as the coarsening continues, i.e., as the discretization with fewer points is considered.   
This phenomenon is known to cause difficulties in designing smoothers to be employed within the MG method for Helmholtz problems~\cite{elman2001multigrid, ernst2013multigrid}.
The most troublesome is the so-called resonance discretization level, in the case of our \nameref{sec:helm1D}  example with $k_{\text{H}}=60$, the level associated with ${\lambda}/{h} \approx 5$.

\begin{figure}
\centering
\includegraphics{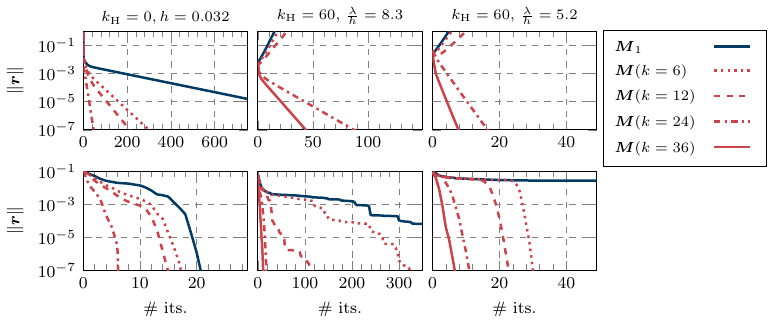}
\caption{The average convergence history of Jacobi preconditioner ($\Mm_1, \gamma_1=\gamma_{k_{\text{H}}} $) and Jacobi preconditioner  multiplicatively composed with coarse space induced by the DeepONet (TB), denoted by~$\Mm$, where $k \in \{6, 12, 24, 36 \}$. 
Preconditioners are used as standalone solvers (top) and for preconditioning of the GMRES(50) (bottom).}
\label{fig:combined_2grid_multiplicative}
\end{figure}

Now, we demonstrate the convergence of the multiplicative preconditioner for all three considered \nameref{sec:helm1D}  test cases. 
In the top row of \cref{fig:combined_2grid_multiplicative}, we illustrate the convergence of our multiplicatively composed preconditioner as a standalone solver. 
For $k_{\text{H}}=0$, augmenting Jacobi smoother with TB-based coarse space increases the convergence speed of the Krylov method.
Notably, larger value of $k$, leads to larger speedup.
In contrast, for $k_{\text{H}}=60$, utilizing a sufficiently large number of TB functions is necessary to ensure the convergence.
This is as expected since we have seen above that the Jacobi iteration diverges for a significant portion of the low-frequency part of the spectrum.
The bottom row of \cref{fig:combined_2grid_multiplicative} demonstrates the convergence history of preconditioned GMRES(50). 
In this case, all considered preconditioners give rise to convergent iteration, independently of the value of~$k$.  
Moreover, compared to a simple, non-hybrid, Jacobi preconditioner, the DeepONet-based Jacobi composite preconditioner is always more efficient.

\subsubsection{Composing level-dependent smoothers for multigrid method}
The two-level preconditioners examined in previous sections are not scalable, especially as the problem size grows. 
To achieve scalability, i.e., the convergence independent of problem size, the MG methods can be employed. 
Although the MG is known to be the optimal preconditioner for elliptic PDEs, it is more challenging to devise the mesh- and level-independent algorithm for different classes of problems, e.g., Helmholtz problems.
We can use the observations made in the previous section to guide the design of the MG, in particular by selecting the level-dependent smoothers. 
To demonstrate such capability, let us consider the \nameref{sec:helm1D}  problem $k_{\text{H}}=60$ discretized using a mesh with $n = 385$ points.
As we have seen in the previous section, the Jacobi smoother tends to amplify the low-frequency components of the error, especially on and around the resonance level. 
As a consequence, it is common practice to employ different smoothers on different levels of the MG. 
The particular choice, involving the type of smoother and number of smoothing steps, typically depends on the ratio~$\lambda/h$. 
A few iterations of the damped Jacobi method are frequently used on all levels, where it converges well. At the same time, the higher-order or Krylov smoothers with a significant number of iterations~\cite{elman2001multigrid} are employed on and around the resonance level. 

{
\begin{table}
\centering
\footnotesize
\begin{tabular}{|c||cccccc||c|}
\hline
\# \textbf{Levels}  & \multicolumn{6}{c||}{\textbf{Smoothing schedule}} & $\#$ \textbf{GMRES its.} \\ \hline
\textbf{2} & 		J & D	& 	& &	&	& \textbf{6}	\\ \hline \hline
\textbf{3} & 		J &J & D	&	& &	&	\textbf{10}	\\ \hline \hline
\textbf{4} & 		J&J &J & D	& & & 	15	\\ \hline
\textbf{4} & 		J &J &  $\Mm (6)$ & D	&	& &   16	\\ 
\textbf{4} & 		J &J & $\Mm (12)$ & D	&	&  & 12	\\ 
\textbf{4} & 		J& J& $\Mm (24)$ & D	&	&  &  \textbf{10}	\\ \hline \hline
\textbf{5} & 		J & J & J & J & D		&	& 28	\\  \hline
\textbf{5} & 		J & J &  J &  $\Mm (6)$ &  D		& & 32	\\
\textbf{5} & 		J & J &  J &  $\Mm (12)$ &   D	& & 16	\\
\textbf{5} & 		J & J &  J &  $\Mm (24)$ & D	& & 15	\\ \hline
\textbf{5} & 		J & J &  $\Mm (6)$ &  $\Mm (6)$ &  D	& & 23	\\
\textbf{5} & 		J & J &  $\Mm (12)$ &  $\Mm (12)$ &  D	& & 14	\\
\textbf{5} & 		J & J &  $\Mm (24)$ &  $\Mm (24)$ & D	& & \textbf{11} \\ \hline \hline
\textbf{6} & 		J & J &  J & J & J	& D	&	28	\\  \hline
\textbf{6} & 		J & J &  J &  $\Mm (6)$ & 		J	& D	& 32	\\
\textbf{6} & 		J & J &  J &  $\Mm (12)$ &   	J 	& D	& 16	\\
\textbf{6} & 		J & J &  J &  $\Mm (24)$ & 	J	& D 	& 15	\\ \hline
\textbf{6} & 		J &J &  $\Mm (6)$ &  $\Mm (6)$ 	&  J 	& D	& 17	\\
\textbf{6} & 		J &J &  $\Mm (12)$ &  $\Mm (12)$ 	&  J 	& D	& 15	\\
\textbf{6} & 		J & J &  $\Mm (24)$ &  $\Mm (24)$ 	& J 	& D	& \textbf{11}	\\ \hline 
\end{tabular}
\caption{The average convergence of GMRES(50) preconditioned with MG for \nameref{sec:helm1D} ($k_{\text{H}}=60$, and $h = 1/384$).
The MG is configured as a V-cycle with varying number of levels and smoothing schedules.
Symbol J denotes $3$ steps of Jacobi with level-dependent $\gamma_{k_{\text{H}}}$, while symbol D stands for direct solver. 
Symbol $\Mm(k)$ denotes composite smoother obtained by multiplicatively composing $3$ steps of Jacobi with DeepONet based coarse space obtained using $k$ randomly selected TB functions.}
\label{tab:MG_k60}
\end{table}}

Here, we use an MG preconditioner with smoothing schedules specified using the following notation:  $\langle J, J, D \rangle$ stands for three-level MG with Jacobi smoother (with $\gamma_{k_{\text{H}}}$) on levels three and two and direct solver (LU factorization) on level one. 
\Cref{tab:MG_k60} summarizes the obtained results. 
If we consider smoothing schedules involving only J and D, the number of GMRES iterations increases as soon as MG has five levels.
This is because the resonance discretization level is now part of the multilevel hierarchy, in particular as a level four. 
As we have seen in the previous section, the Jacobi smoother amplifies the low-frequency components of the error on this level, which in turn hinders the overall convergence of the MG method.

To enable the level-independent convergence of the MG method, we augment Jacobi with a DeepONet.
Following the discussion from the previous section, we do so by multiplicatively composing Jacobi with a TB-based coarse space operator. 
Since such composite smoother is computationally more expensive than simple Jacobi, we use it only around the resonance level.
The results reported in \cref{tab:MG_k60}, demonstrate that incorporating composite smoother gives rise to level-independent MG, provided that TB-based coarse space is sufficiently big. 
For example, let us consider the MG with five levels. 
If we augment the second coarsest level with the DeepONet component, we obtain convergence comparable to the MG method with four levels (Jacobi smoothers), i.e., the GMRES method converges with $15$ iterations. 
This is expected, as the composite smoother effectively reduces only the error components associated with this particular level. 
If we further augment the third coarsest level, we can now achieve convergence comparable to the three-level method (Jacobi smoothers), i.e., the GMRES method convergence with approximately $10$ iterations. 
Thus, we have achieved level-independent convergence by mitigating the drawbacks of Jacobi smoother around the resonance level.

The quality of the proposed smoothers can be further improved by carefully selecting the ratio between a number of Jacobi and DeepONet-based steps and by finding optimal scaling parameters~$\gamma_1, \gamma_2$. 
Moreover, for the practicality of the proposed approach, it is crucial to devise adaptive strategies for automatically selecting the value of $k$ on each level. 
The results obtained in this section demonstrate a clear potential for exploring such research directions in the future.

 \subsection{Additive preconditioning}
\label{sec:results_additive}
We continue our numerical study by investigating the performance of the two-level ASM preconditioner as specified by~\eqref{eq:two_level_additive}, where we set $\gamma_s=1$, for all $s=0, \ldots, S$. 
The subdomains are constructed using the METIS partitioner~\cite{karypis1997metis}.
The coarse space is created using the DeepONet. 

\begin{table}
\centering
\footnotesize
\begin{tabular}{|l|c||r|r|r|r|r|r|r|r|r|r|r|r|r|r|r|r|r|}
\hline
\multicolumn{2}{|c||}{ $\mathbf{S}$}      & \multicolumn{2}{c|}{$\mathbf{16}$}                &\multicolumn{2}{c|}{$\mathbf{32}$}                  & \multicolumn{2}{c|}{$\mathbf{64}$}         \\ \cline{1-8}
\multicolumn{2}{|c||}{ $\mathbf{k /N_S}$}	& $\mathbf{1,000}$   & $\mathbf{2,500}$   & $\mathbf{1,000}$  & $\mathbf{2,500}$  &  $\mathbf{1,000}$    & $\mathbf{2,500}$  \\ \hline 
 \multirow{4}{*}{\rotatebox[origin=c]{90}{$\pazocal{T}^1$-dense}} & {0} 
& \multicolumn{2}{c|}{$60.2 \pm 0.6$} 	 
& \multicolumn{2}{c|}{$82.1 \pm 1.4$} 	
& \multicolumn{2}{c|}{$84.9 \pm 1.3$}       \\ \cline{2-8}
&{8} 
& $42.4 \pm 0.8$	& $42.8 \pm 1.5$
& $57.6 \pm 1.7$ 	& $59.5 \pm 1.7$ 
&$57.3 \pm 2.6$   & $57.3 \pm 58.8$     \\ \cline{2-8}
 &{32} 
& $26.9 \pm 0.3$ 	&	$25.9 \pm 0.3$   
& $35.4 \pm 0.5$  	& 	$34.2 \pm 0.6$  
& $35.9 \pm 0.7$ 	& 	$33.4 \pm 0.5$  \\ \cline{2-8}
&{128}
&	$19.3 \pm 0.5$ 		& 	$17.1 \pm 0.3$  
& 	$20.8 \pm 0.4$  	& 	$18.2 \pm 0.4$
&  	$21.5 \pm 0.7$  	& 	$18.3 \pm 0.45$  \\ \hline  \hline
 \multirow{4}{*}{\rotatebox[origin=c]{90}{$\pazocal{T}^1$-sparse}}
&{1} 
& 	$53.4 \pm 1.0$ 	& 	$53.8 \pm 0.7$   
&	$78.2 \pm 0.9$ 	& 	$78.0 \pm 1.3$ 
& 	$71.0 \pm 1.2$   & 	$71.2 \pm 1.9$     \\ \cline{2-8}
&{3} 
& $44.8 \pm 1.4$ 	& $45.8 \pm 2.4$   
& $60.0 \pm 0.9$  	& $65.0 \pm 3.5$  
& $40.6 \pm 1.0$ 	& $43.0 \pm 0.6$  \\ \cline{2-8}
&{5}
& $37.8 \pm 1.3$ 	& $36.4 \pm 1.3$   
& $43.4 \pm 1.9$  	& $46.0 \pm 2.5$  
& $26.6 \pm 0.8$ 	& $26.2 \pm 0.7$  \\ \cline{2-8}
&{8}
& $31.0 \pm 0.1$ 	& $28.4 \pm 0.5$   
& $31.8 \pm 0.7$  	& $29.4 \pm 0.5$  
&  $17.8 \pm 0.8$  & $17.2 \pm 0.4$  \\ \hline \hline
\multicolumn{2}{|c||}{ $\mathbf{S}$}      & \multicolumn{2}{c|}{$\mathbf{16}$}                &\multicolumn{2}{c|}{$\mathbf{32}$}                  & \multicolumn{2}{c|}{$\mathbf{64}$}         \\ \cline{1-8}
\multicolumn{2}{|c||}{ $\mathbf{k /N_S}$}	& $\mathbf{1,000}$   & $\mathbf{2,500}$   & $\mathbf{1,000}$  & $\mathbf{2,500}$ &  $\mathbf{1,000}$    & $\mathbf{2,500}$   \\ \hline 
 \multirow{4}{*}{\rotatebox[origin=c]{90}{$\pazocal{T}^4$-dense}}&{0} 
& \multicolumn{2}{c|}{$262.3 \pm 1.7$} 	 
& \multicolumn{2}{c|}{$294.4 \pm 3.4$} 	
& \multicolumn{2}{c|}{$366.1 \pm 3.5$}       \\ \cline{2-8}
&{8} 
& $181.9 \pm 5.2$ 	& $183.6 \pm 4.7$   
& $211.5 \pm 7.3$ 	& $208.4 \pm 4.9$ 
& $257.4 \pm 6.3$   & $256.2 \pm 8.6$     \\ \cline{2-8}
& {32} 
& $112.8 \pm 1.1$ 	& $107.5 \pm 1.6$   
& $127.8 \pm 1.7$  	& $121.2 \pm 1.8$  
& $157.2 \pm 1.5$ 	& $147.9 \pm 1.6$  \\ \cline{2-8}
&{128 }
& $\phantom{0}60.1 \pm 0.6$ 	& $\phantom{0}56.7 \pm 0.9$   
& $\phantom{0}64.1 \pm 0.8$  	& $\phantom{0}63.2 \pm 0.6$  
&  $\phantom{0}71.2 \pm 1.2$ 	& $\phantom{0}72.6 \pm 0.9$  \\ \hline  \hline
 \multirow{4}{*}{\rotatebox[origin=c]{90}{$\pazocal{T}^4$-sparse}}
&{1} 
& 	$213.8 \pm 2.3$ 	& 	$213.8 \pm 1.2$   
& 	$223.7 \pm 7.3$ 	& 	$224.8 \pm 5.1$ 
&	$245.5 \pm 6.3$   & 		$245.5 \pm 5.9$     \\ \cline{2-8}
&{3} 
&$155.4 \pm 3.3$ 	& $155.2 \pm 3.1$   
& $160.4 \pm 6.2$ 	& $160.0 \pm 5.8$ 
& $174.2 \pm 5.8$   & $174.2 \pm 4.2$     \\ \cline{2-8}
& {5} 
&$134.6 \pm 1.3$ 	& $134.4 \pm 3.1$   
& $135.6 \pm 1.7$  	& $135.7 \pm 1.2$  
& $144.5 \pm 1.7$ 	& $144.5 \pm 1.3$  \\ \cline{2-8}
&{8}
& $107.5 \pm 0.2$ 	& $107.0 \pm 0.3$   
& $109.2 \pm 0.5$  	& $109.2 \pm 0.5$  
& $119.6 \pm 0.7$  	& $119.0 \pm 0.4$  \\ \hline 
\end{tabular}
\caption{The average number of iterations required by the CG, preconditioned using ASM ($S$ subdomains) with additive coarse space.
The coarse space is constructed using the TB approach with $k$ TB functions, where $k=0$ indicates no coarse space.
The experiments performed for  \nameref{sec:diff} example, discretized using $\pazocal{T}^{1}$ and $\pazocal{T}^{4}$ meshes.}
\label{tab:additive_ASM_summary}
\end{table}

Our first set of experiments in this section considers the \nameref{sec:diff} example.
Firstly, we investigate the performance of the ASM preconditioner with the coarse space constructed using the TB approach, with the dense and sparse transfer operators, as discussed in \cref{sec:sparsity_dd}. 
In the case of the sparse transfer operators, one step of Jacobi prolongation smoothing is utilized. 
\Cref{tab:additive_ASM_summary} reports the obtained results with respect to an increasing number of subdomains for two different problem sizes related to meshes~$\pazocal{T}^1$ and~$\pazocal{T}^4$.
As we can see from the observed results, increasing the number of samples does not significantly improve the results, suggesting that the DeepONets can be trained with few samples, i.e., at low computational cost. 
We also see that using the dense approach requires a larger number of TB functions to achieve algorithmic scalability, i.e., convergence with a fixed number of iterations for an increasing number of subdomains. 
For instance, in the case of a dense approach, it is necessary to utilize $128$ TB functions for a number of iterations to remain bounded with respect to the increasing number of subdomains. 
In contrast, utilizing a sparse approach, it suffice to use fewer TB functions to obtain an algorithmically scalable method.
Moreover, we point out that our numerical results demonstrate that the DeepONet-based two-level ASM is algorithmically scalable even for the problem discretized using $\pazocal{T}^4$. 
This highlights the generalization capabilities of the proposed preconditioning framework, as the DeepOnet has been trained using the data generated with $\pazocal{T}^1$.

We also compare the performance of the two-level ASM preconditioner  with coarse space obtained using TB (sparse, $k=8$) and DP approaches. 
\Cref{tab:additive_ASM_summary_DP} summarizes the obtained results. 
Here, we recall that although the operator~$\Am$ is SPD, the DP approach requires us to employ the F-GMRES method. 
The results are reported with respect to the increasing number of samples and with an increasing number of subdomains. 
To use the DP approach, we sample the right-hand side using random distribution~$\sim \pazocal{U}[0, 1]$.
The obtained results are in accordance with the observations made in the previous section. 
As we can see, using a standalone ASM preconditioner, the F-GMRES(50) method requires more iterations to achieve convergence than the CG method due to the restarting procedure. 
In contrast to the TB approach, the performance and algorithmic scalability of the DP-based preconditioner improve drastically by incorporating a larger number of samples. 
These results are consistent with findings reported in \cref{sec:results_multiplicative} and suggest that investing more resources into sampling and training procedures gives rise to a more robust DP-based preconditioner.

{\begin{table}
\centering
\footnotesize
\begin{tabular}{|l||l||r|r|r|r|r|r|r|}
\hline
\multicolumn{2}{|c||}{ $\mathbf{S}$ }    & \multicolumn{3}{c|}{$\mathbf{16}$}             \\ \hline
\multicolumn{2}{|c||}{$\mathbf{N_S}$}    & $\mathbf{2,500}$   & $\mathbf{25,000}$   & $\mathbf{100,000}$  \\ \hline
 \multirow{4}{*}{${\bpazocal{T}^1}$}  
 & CG (ASM) & $60.2 \pm 0.6\phantom{0}$ & $60.2 \pm 0.6\phantom{0}$ & $60.2 \pm 0.6\phantom{0}$  \\ \cline{2-5}
& CG (ASM+TB)
& $28.4 \pm 0.5\phantom{0}$ & $28.2 \pm 0.2\phantom{0}$  & $28.7 \pm 0.1\phantom{0}$  \\ \cline{2-5}  
 & F-GMRES (ASM) & $89.2 \pm 3.2\phantom{0}$  & $89.2 \pm 3.2\phantom{0}$  & $89.2 \pm 3.2\phantom{0}$  \\ \cline{2-5}
& F-GMRES (ASM+DP) & $108.2 \pm 13.2$ & $57.8 \pm 6.9\phantom{0}$  & $41.2 \pm 5.3\phantom{0}$  \\ \hline  \hline
 \multirow{4}{*}{$\bpazocal{T}^4$}  
 & CG (ASM) &  $262.3 \pm 1.7\phantom{0}$  &  $262.3 \pm 1.7\phantom{0}$  &  $262.3 \pm 1.7\phantom{0}$  \\ \cline{2-5}
&  CG (ASM+TB) & $107.0 \pm 0.3\phantom{0}$ & $102.2 \pm 0.2\phantom{0}$  & $105.7 \pm 0.3\phantom{0}$ \\ \cline{2-5} 
&  F-GMRES (ASM) &  $391.2 \pm 3.5\phantom{0}$  &  $391.2 \pm 3.5\phantom{0}$ &  $391.2 \pm 3.5\phantom{0}$  \\ \cline{2-5}
& F-GMRES (ASM+DP) & $452.1 \pm 23.4$ & $353.1 \pm 14.8$  & $245.1 \pm 10.2$  \\ \hline  \hline
\multicolumn{2}{|c||}{ $\mathbf{S}$ }    & \multicolumn{3}{c|}{$\mathbf{32}$}             \\ \hline
\multicolumn{2}{|c||}{$\mathbf{N_S}$}    & $\mathbf{2,500}$   & $\mathbf{25,000}$   & $\mathbf{100,000}$  \\ \hline
 \multirow{4}{*}{$\bpazocal{T}^1$}  
 &  CG (ASM) &   $82.1 \pm 1.4\phantom{0}$  &   $82.1 \pm 1.4\phantom{0}$  &   $82.1 \pm 1.4\phantom{0}$  \\ \cline{2-5}
 & CG (ASM+TB) &  $29.4 \pm 0.5\phantom{0}$ &  $29.5 \pm 0.1\phantom{0}$ &  $29.3 \pm 0.2\phantom{0}$  \\ \cline{2-5} 
 & F-GMRES (ASM) &   $112.2 \pm 2.3\phantom{0}$ &   $112.2 \pm 2.3\phantom{0}$ &   $112.2 \pm 2.3\phantom{0}$ \\ \cline{2-5}
& F-GMRES (ASM+DP) 	&  $134.2 \pm 10.2$ &  $104.2 \pm 5.5\phantom{0}$ &  $62.1 \pm 2.1\phantom{0}$  \\ \hline  \hline
 \multirow{4}{*}{$\bpazocal{T}^4$}  
 & CG (ASM)  & $294.4 \pm 3.4\phantom{0}$ & $294.4 \pm 3.4\phantom{0}$ & $294.4 \pm 3.4\phantom{0}$ \\ \cline{2-5}
& CG (ASM + TB) 		&  $109.2 \pm 0.5\phantom{0}$ &  $108.4 \pm 0.3\phantom{0}$ &  $109.5 \pm 0.2\phantom{0}$ \\ \cline{2-5} 
 & F-GMRES (ASM)  	& $492.3 \pm 9.2\phantom{0}$ & $492.3 \pm 9.2\phantom{0}$  & $492.3 \pm 9.2\phantom{0}$  \\ \cline{2-5}
& F-GMRES (ASM+DP) 	&  $569.9 \pm 23.1$ &  $479.7 \pm 12.4$ &  $286.3 \pm 9.9\phantom{0}$  \\ \hline \hline 
\multicolumn{2}{|c||}{ $\mathbf{S}$ }    & \multicolumn{3}{c|}{$\mathbf{64}$}             \\ \hline
\multicolumn{2}{|c||}{$\mathbf{N_S}$}    & $\mathbf{2,500}$   & $\mathbf{25,000}$   & $\mathbf{100,000}$  \\ \hline
 \multirow{4}{*}{$\bpazocal{T}^1$}  
 & CG (ASM) &   $84.9 \pm 1.3\phantom{0}$  &   $84.9 \pm 1.3\phantom{0}$  &   $84.9 \pm 1.3\phantom{0}$  \\ \cline{2-5}
& CG (ASM+TB)    &  $\phantom{0}17.8 \pm 0.8\phantom{0}$ &  $17.2 \pm 0.4\phantom{0}$ &  $17.3 \pm 0.5\phantom{0}$  \\ \cline{2-5} 
 & F-GMRES (ASM) &   $125.2 \pm 2.7\phantom{0}$ &   $125.2 \pm 2.7\phantom{0}$ &   $125.2 \pm 2.7\phantom{0}$ \\ \cline{2-5}
& F-GMRES (ASM+DP) &  $123.4 \pm 8.1\phantom{0}$ &  $111.2 \pm 4.9\phantom{0}$ &  $85.3 \pm 4.1\phantom{0}$  \\ \hline  \hline
 \multirow{4}{*}{$\bpazocal{T}^4$}  
 & CG (ASM)  & $366.1 \pm 3.5\phantom{0}$ & $366.1 \pm 3.5\phantom{0}$ & $366.1 \pm 3.5\phantom{0}$ \\ \cline{2-5}
& CG (ASM+TB) &  $119.6 \pm 0.7\phantom{0}$ &  $119.0 \pm 0.4\phantom{0}$ &  $119.3 \pm 0.3\phantom{0}$  \\ \cline{2-5} 
 & F-GMRES (ASM)  & $521.3 \pm 10.6$ & $521.3 \pm 10.6$  & $521.3 \pm 10.6$  \\ \cline{2-5}
& F-GMRES (ASM+DP) &  $587.4 \pm 25.2$ &  $475.2 \pm 15.2$ &  \multicolumn{1}{c|}{$287.4 \pm 10.2$}  \\ \hline 
\end{tabular}
\caption{The average number of iterations required by the Krylov methods preconditioned with ASM ($S$ subdomains) to reach convergence. 
The additive DeepONet-based coarse space is incorporated using TB (sparse, $k=8$) and DP (right-hand side sampled from~$\sim \pazocal{U} [0, 1]$) approaches.
The DeepONets are trained using~$N_S$ samples.
The experiments performed for \nameref{sec:diff} example discretized using $\pazocal{T}^{1}$ and $\pazocal{T}^{4}$ meshes.}
\label{tab:additive_ASM_summary_DP}
\end{table}}

\begin{figure}
  \centering
\includegraphics{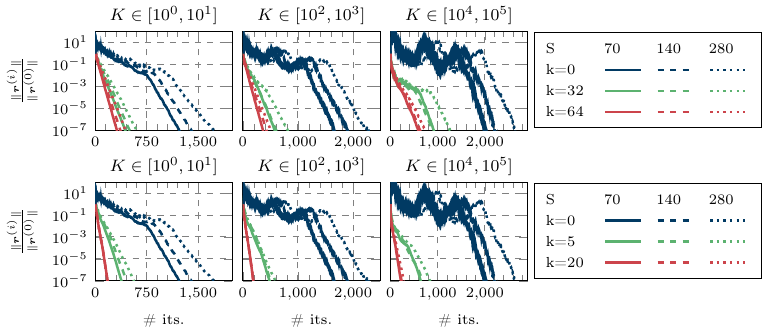}  
\caption{The average convergence of the CG algorithm, preconditioned using ASM with the overlap of two for the \nameref{sec:jump_diff} example, discretized using $\pazocal{T}^6$. 
The number of ASM's subdomains (S) is selected from $\{70, 140, 280 \}$.
The coarse space is constructed using dense $\Pm$ (top) and sparse $\Pm$ (bottom)  obtained by extracting $k$ TB functions.}
\label{fig:jump_asm_adaptive}
\end{figure}

Our second example considers the \nameref{sec:jump_diff} problem. 
\Cref{fig:jump_asm_adaptive} demonstrates the obtained results for the problems discretized using $\pazocal{T}^6$ mesh with an increasingly larger jump in the diffusion coefficient.  
The top row of \cref{fig:jump_asm_adaptive} depicts results obtained using dense transfer operators, constructed with $32$ and $64$ TB functions. 
The experiments are performed for an increasing number of subdomains to study the scalability of the resulting preconditioner. 
As we can see, using a larger number of TB functions improves the method's convergence and scalability. 
The bottom row of \cref{fig:jump_asm_adaptive} depicts results obtained using sparse, transfer operators constructed with $5$ and $20$ TB functions. 
The obtained results demonstrate that using $20$ TB vectors is sufficient to obtain an efficient and algorithmically scalable algorithm.

\begin{figure}
  \centering
\includegraphics{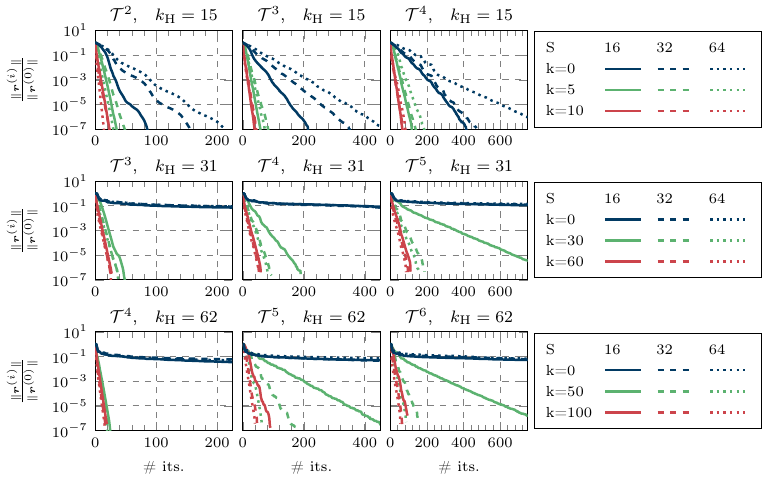}    
\caption{The average convergence of the GMRES(50), preconditioned using ASM without overlap for the parametric \nameref{sec:helm2D} problem of varying sizes and $k_{\text{H}}$ numbers.
The number of ASM's subdomains (S) is selected from the set $\{16, 32, 64 \}$. 
The DeepONets are trained using $\pazocal{T}^{1-3}$ meshes for $k_{\text{H}} \in \{15, 31, 62 \}$, respectively.
The coarse space is constructed using sparse $\Pm$, obtained by utilizing TB functions.}
\label{fig:helm_asm_adaptive}
\end{figure}

Ultimately, we test our two-level DeepONet-based ASM preconditioner for \nameref{sec:helm2D} problem with $k_{\text{H}} \in \{ 15, 31, 62 \}$. 
\Cref{fig:helm_asm_adaptive} demonstrates the results with respect to increasing resolution of the problem.
We can see that using a single-level ASM preconditioner does not give rise to a convergent method for $k_{\text{H}}=31$ and $k_{\text{H}}=62$. 
The divergence is concluded after the method fails to converge within $8,000$ iterations. 
The results also demonstrate that as $k_{\text{H}}$ increases, a larger number of TB functions is required to obtain an algorithmically scalable two-level ASM.
For example, for $k_{\text{H}}=15$, we require around $10$ TB functions, while for $k_{\text{H}}=31$, it is necessary to utilize around $60$ TB functions.

The numerical results demonstrate that we can construct effective coarse space for ASM by using the TB approach. 
Compared to the standard Nicolaides approach~\cite{dolean2015introduction}, the resulting coarse space is computationally more expensive. 
However, in contrast to the Nicolaides approach, it is effective even in the presence of jumping coefficients or indefiniteness. 
Compared to more advanced coarse space approaches, such as the Dirichlet-to-Neumann map~\cite{dolean2012analysis}, or GenEO~\cite{spillane2014abstract}, the construction of coarse space using the TB approach is significantly simplified and less computationally demanding. 
We also note that the obtained results could be further improved by incorporating the physics-based priors into the trunk network architecture and by designing more elaborate techniques for selecting the TB functions from the DeepONet.

 \section{Conclusion}
\label{sec:conclusion}
We proposed a novel DeepONet-based preconditioning framework for solving parametric linear systems of equations. 
The devised preconditioners were constructed using the subspace correction framework and leveraged the DeepONet and standard iterative methods to eliminate low and high-frequency components of the error, respectively. 
The DeepONet hybridization was performed using DP and TB approaches. 
The DP approach utilized a DeepONet inference to approximate the inverse of the parametrized linear operator at each preconditioning step.
In contrast, the TB approach used DeepONet to construct the subspace problem, the solution of which was obtained using standard numerical methods. 
Performed numerical results demonstrated that the TB approach outperforms the DP approach. 
Furthermore, we have shown that using the TB approach gives rise to coarse space, which is trivial to construct and enables the construction of algorithmically scalable ASM.

We foresee the extension of the presented work in several ways.
For instance, the efficacy of the DP approach could be improved by developing advanced sampling strategies for the right-hand sides encountered during the preconditioning step. 
This improvement might be achieved, for example, by training the DeepONet using a solver-in-loop approach~\cite{um2020solver}. 
Additionally, the TB approach could be enhanced by adopting more elaborate strategies to select an optimal set of basis functions from the DeepONet. 
The quality of these basis functions could be further improved by integrating prior knowledge into the architectural design.
Moreover, innovative aggregation strategies can be devised to enforce the sparsity pattern in the resultant transfer operator. 
It would also be interesting to investigate the convergence properties of the proposed preconditioning framework, considering different types of operator learning approaches, e.g., CNO~\cite{raonic2023convolutional}, and varying geometries.
The applicability of the proposed preconditioning framework could be further extended to different problem types, such as nonlinear or eigenvalue problems.
Moreover, the current software implementation could be expanded to distributed settings.

\appendix
\section{Representative examples of Krylov methods}
\label{sec:krylov_examples}
This section provides algorithmic details for Krylov methods, which we use to generate the numerical results reported in \cref{sec:num_results}, namely GMRES and CG. 

\subsection{Flexible generalized minimal residual method}
We consider a particular variant of GMRES\footnote{The GMRES(m) algorithm can be obtained from \cref{alg:fgmres} by considering constant preconditioner~$\Mm$ and setting~$\Zm^{(i)} := \Vm^{(i)}$. }, called flexible GMRES (F-GMRES)~~\cite{saad1993flexible}, which allows for variable preconditioning. 
In this particular case, the basis of the subspace $\pazocal{K}_i(\Am, \rv^{(0)})$ are constructed using the flexible Arnoldi procedure$\footnote{In case of standard GMRES, $\Zm^{(i)} := \Vm^{(i)}$ is utilized by the Arnoldi procedure.}$, which utilizes the following decomposition:
\begin{align}
\Am \Zm^{(i)} = \Vm^{(i+1)} \widebar{\Hm}^{(i)} \qquad \text{with} \qquad  (\Vm^{(i+1)})^{\top} \Vm^{(i+1)}  = \Id, 
\label{eq:arnoldi}
\end{align} 
where $\widebar{\Hm}^{(i)} \in \R^{(i+1) \times i}$ denotes an upper Hessenberg matrix.
The matrix ${\Zm^{(i)} = [\zv^{(1)}, \ldots, \zv^{(i)}] \in \R^{n \times i}}$ is orthonormal, and its first column consists of normalized residual, i.e.,~$\zv^{(1)}=\rv^{(0)}/\zeta$, where $\zeta = \| \rv^{(0)} \|$. 
Subsequent columns are constructed by orthogonalizing~$\Am \zv^{(i)}$ with respect to previous basis vectors.

After the basis of~$\pazocal{K}_i(\Am, \rv^{(0)})$ are constructed, the solution~$\uv^{(i)}$ can be found by minimizing the residual norm over the subspace~$\uv^{(0)} + \text{range}(\Zm^{(i)})$, i.e., 
\begin{align}
\uv^{(i)}:= \underset{\uv^{(i)} \in \uv^{(0)} + \pazocal{K}_i(\Am, \rv^{(0)})}{\text{arg min}} & \| \fv - \Am \underbrace{(\uv^{(0)} + \Zm^{(i)}\ym^{(i)})}_{\uv^{(i)}}  \|_2, 
\label{eq:functional_gmres}
\end{align}
where~$\ym^{(i)} \in \R^{i}$. 
Here, we point out that the term $\| \fv - \Am \uv^{(i)} \|_2$  can be expressed as ${ \| \zeta \ev_1 -   \widebar{\Hm}^{(i)} \ym^{(i)}    \|_2 }$, where $\ev_1 = [1, 0, \ldots, 0]^{\top} $ is the first canonical basis vector in $\R^i$.

In the end, we note that for the large-scale applications, the restarted variants of F-GMRES~\cite{morgan1995restarted}, denoted F-GMRES(m),  have to be employed in order to reduce the memory footprint associated with storing $\Zm^{(i)}$ and $\Vm^{(i)}$. 
The restarting strategy simply relaunches the algorithm from scratch on every $m$-th iteration, starting from the current solution approximation. 
The F-GMRES(m) algorithm is summarized in \cref{alg:fgmres}.

We note that no general convergence results for the F-GMRES method are available, as the subspace $\pazocal{K}_i$ is no longer a Krylov space, but it is nevertheless the space in which the solution is sought.
However, the analysis of F-GMRES(m) breakdown, similar to that of the GMRES method, can be found, for example, in~\cite{saad1993flexible, saad2003iterative}.

\subsection{Conjugate gradient method}
Using CG method, the basis of the subspace $\pazocal{K}_i(\Am, \rv^{(0)})$ are constructed using Arnoldi like procedure, called Lanczos method~\cite{saad2003iterative}, which utilizes~$\Zm^{(i)} = \Vm^{(i)}$. 
At each iteration, the approximate solution $\uv^{(i)}$ is computed as $\uv^{(i)} = \uv^{(0)} + \Vm^{(i)} \ym$. 
Here, the vector $\ym$ is obtained by solving the small linear system~$\Hm^{(i)} \yv = \zeta \ev_1$, where ${\Hm^{(i)} = (\Vm^{(i)})^{\top}  \Am \Vm^{(i)}}$.
Since~$\Am$ is SPD, the matrix~$\Hm^{(i)}$ is tridiagonal.
As a consequence, the orthogonalization of~$\Vm^{(i)}$ can be performed cheaply using three-term recurrence. 
Moreover, a LU factorization of~$\Hm^{(i)}$ can be computed incrementally, and it is known to be numerically stable without pivoting.

Taking advantage of these properties, the CG algorithm can be implemented using short recurrences that do not require storage of~$\Vm^{(i)}$, see \cref{alg:deflated_pcg}. 
In particular, at each iteration, the orthogonality and conjugacy conditions are utilized to construct the approximation~$\uv^{(i)}$ as follows
\begin{align}
\uv^{(i)} = \uv^{(i-1)} + \alpha^{(i-1)} \pv^{(i-1)}.
\end{align}
The search directions are constructed recursively. 
In particular~$\pv^{(i)}$ is given by a linear combination of~$\rv^{(i)}$ and the previous search direction~$\pv^{(i-1)}$, i.e., 
\begin{align}
\pv^{(i)} = \rv^{(i)} + \beta^{(i-1)} \pv^{(i-1)},   \quad \text{for}  \ i \geq 1,
\end{align}
and~$\pv^{(i)} = \rv^{(i)}$, for $i=0$.
The parameter~$\alpha^{(i-1)}$ is chosen as  ${\alpha^{(i-1)}=\frac{\langle  \rv^{(i-1)},  \rv^{(i-1)} \rangle}{\langle \pv^{(i-1)},  \Am \pv^{(i-1)} \rangle}}$, i.e., such that the residuals~$\rv^{(i)}$ and $\rv^{(i-1)}$ are orthogonal to each other.  
In addition, the parameter~$\beta^{(i-1)}$ is obtained by enforcing the conjugacy between~$\pv^{(i)}$ and~$ \pv^{(i-1)}$, thus as~${\beta^{(i-1)}= \frac{\langle  \rv^{(i)},   \rv^{(i)} \rangle}{\langle  \rv^{(i-1)},  \rv^{(i-1)} \rangle}}$.

The CG algorithm is well-known for its computational efficiency and low memory requirements.
Moreover, an upper bound on the convergence rate exists, which states an explicit dependence on a condition number of~$\Am$, see~\cite{meurant2006lanczos, saad2003iterative}.

\begin{algorithm}[t]
  \caption{Preconditioned Restarted Flexible GMRES (F-GMRES(m))}
 \label{alg:fgmres}
  \begin{algorithmic}[1]
\Require {$\Am \in \R^{n \times n}, \fv, \uv^{(0)} \in \R^n, \Mm \in \R^{n \times n}$, $m \in \N$, $\widebar{\Hm} \in \R^{m \times m}$}
\While{not converged} 
\State $\zeta=\| \rv^{(0)} \|$, $\vv^{(1)} = \rv^{(0)}/\zeta$, where $\rv^{(0)}  = \fv - \Am \uv^{(0)}$
\State $\widebar{\Hm} = \boldsymbol{0}$
\For{i=1, \ldots,  m}     \COMMENTmine{Arnoldi process}
\State $\zv^{(i)} = \Mm \vv^{(i)}$	 \COMMENTmine{Flexible preconditioning step}
\State $\wv = \Am \zv^{(i)}$
\For{j=1, \ldots,  i}    
\State $\widebar{\Hm}_{j,i} = \langle \wv , \vv^{(j)} \rangle$
\State $\wv = \wv - \widebar{\Hm}_{j,i}  \vv^{(j)}$
\EndFor
\State $\widebar{\Hm}_{i+1, i} = \| \wv \| $ \ and  \ $\vv^{(i+1)} = \wv/\| \wv \|$
\EndFor
\State $\Zm^{(m)} = [ \zv^{(1)}, \ldots, \zv^{(m)} ]$
\State
\State $\ym^{(m)} = \text{argmin}_{\yv \in \R^{n}} \ \| \zeta \ev_1 - \widebar{\Hm} \yv \|_2  $  \COMMENTmine{Update solution}
\State $\uv^{(m)} = \uv^{(0)} +\Zm^{(m)} \yv^{(m)}$
\State $\uv^{(0)} = \uv^{(m)}$ \COMMENTmine{Restart}
\EndWhile
\State
\Return  $\uv^{(m)}$
\end{algorithmic}
\end{algorithm}

\begin{algorithm}[t]
  \caption{Preconditioned Conjugate Gradient (PCG)}
 \label{alg:deflated_pcg}
  \begin{algorithmic}[1]
\Require {$\Am \in \R^{n \times n}, \fv \in \R^n, \uv^{(0)} \in \R^n, \Mm \in \R^{n \times n}, \Rm \in \R^{k \times n}, $}
\State $\zv^{(0)} = \Mm \rv^{(0)} $, where $\rv^{(0)}  = \fv - \Am \uv^{(0)}$
\State $\pv^{(0)} = \zv^{(0)}$ 
\While{$i=1, 2, \ldots,$ \textbf{until} convergence} 
\State $\alpha^{(i-1)} = \langle \rv^{(i-1)}, \zv^{(i-1)} \rangle/ \langle \pv^{(i-1)},  \Am \pv^{(i-1)} \rangle$
\State $\uv^{(i)} = \uv^{(i-1)} + \alpha^{(i-1)}  \pv^{(i-1)}$
\State $\rv^{(i)} = \rv^{(i-1)} - \alpha^{(i-1)}  \Am \pv^{(i-1)}$
\State $\zv^{(i)} = \Mm \rv^{(i)} $   \COMMENTmine{Preconditioning step}
\State $\beta^{(i-1)} =  \langle \rv^{(i)}, \zv^{(i)} \rangle/ \langle \rv^{(i-1)},  \zv^{(i-1)} \rangle$
\State$ \pv^{(i)} = \beta^{(i-1)}  \pv^{(i-1)} + \zv^{(i)}$ 
\EndWhile
\State
\Return  $\uv^{(i)}$
\end{algorithmic}
\end{algorithm}

\section{Restriction operator for DP approach}
\label{sec:appendix_mapping}
In this section, we discuss how to assemble transfer operator~$\Rm$, utilized by the HINTS/DP approach discussed in~\eqref{eq:q_dd_2l_a}. 
Let $r$ be the index of the branch network, which encodes the parametrized right-hand side. 
The transfer operator~$\Rm$ has to be designed such that it maps the residual~$\rv^{(i)}$ from the dual FE space $\pazocal{V}_h'$ to the space~$\pazocal{Y}^r$.
This is due to the fact that for a given model parameters~$\boldsymbol{\theta}$ and the current nodal approximation~$\uv^{(i)}$, the residual $\rv^{(i)} \in \R^n$ is given as
\begin{align*}
\rv^{(i)} 	= \fv  - \Am {\uv^{(i)}} = (f(v_h; \boldsymbol{\theta}) -   a({u_h^{(i)}}(\boldsymbol{\theta}) , v_h; \boldsymbol{\theta}) ) 	= \langle {r}^{(i)}_h({u_h^{(i)}}; \boldsymbol{\theta}), v_h \rangle,
\end{align*}
for all~$v_h \in \pazocal{V}_h$. 
Here, ${r}^{(i)}_h$ is evaluated as  ${r}^{(i)}_h = \sum_{j=1}^n \psi_j(x) \rv_{\vv}^{(i)}$, where $\rv_{\vv}^{(i)}$ contains the nodal values of the residual, for a given~$\uv^{(i)}$.

As a consequence, $\rv^{(i)}$ is given as 
$${\rv^{(i)} = \bigg{\langle} \sum_{k=1}^n \psi_k(x) (\rv_{\vv}^{(i)})_k, \sum_{j=1}^n \psi_j(x) \bigg{\rangle}}.$$ 
This can be algebraically expressed as $\rv^{(i)} = \Mm_{\mm} \rv_{\vv}^{(i)}$, 
where $\Mm_{\mm} \in \R^{n \times n}$ denotes the mass matrix, elements of which are defined as~$(\Mm_{\mm})_{kj} = \langle \psi_k, \psi_j \rangle$. 
Thus, for a given~$\rv^{(i)}$, the nodal coefficients~$\rv_{\vv}^{(i)}$ can be retrieved as $\rv_{\vv}^{(i)} = \Mm_{\mm}^{-1} \rv^{(i)}$.
Although this operation requires an inverse of the mass matrix, it can be performed efficiently since the mass-matrix $\Mm_{\mm}$ is typically well-conditioned. 
Moreover, in case of the first-order FE, we can approximate $\Mm_{\mm}$ by lumped mass matrix, which is diagonal and therefore trivial to invert. 
Once, the nodal values~$\rv_{\vv}$ are obtained, they can be mapped to $\pazocal{Y}^r$ spanned by the basis functions~$\{ \phi_j \}_{j=1}^{nb_r}$. 
This can be, for example, achieved by using (higher-order) interpolation strategies.

We note that this transformation, involving the mass-matrix of a high-fidelity problem, is only required if the PDE and branch input features are discretized using non-uniform meshes or if the basis functions~$\{ \phi_j \}_{j=1}^{nb_r}$ and $\{ \psi \}_{j=1}^{n}$ are not nodal. 
If the meshes are uniform, and the basis functions of~$\pazocal{V}_h$ and~$\pazocal{Y}^r$ are nodal, then the application of $\Mm_{\mm}^{-1}$ provides a uniform scaling, which the normalization of the branch input can absorb.

\section{Numerical approximation details}
\label{sec:appB}
In this section, we provide details associated with high and low-fidelity numerical approximations of~\eqref{eq:param_pde}. 
In particular, \cref{table:mesh_hierarchies} gives a summary information regarding the meshes used for FE discretization. 
Unless specified differently, we use meshes~$\pazocal{T}^1$ to construct the dataset required for training DeepONets. 
The details of the DeepONets' architectures are shown in \cref{tab:architecture}. 
For convolutional networks (Conv), we always choose the kernel size to be three, while stride is set to two.
The feed-forward networks (FFNs) employ standard dense layers consisting of weights and biases.
\Cref{table:training_params} displays the training times for DeepONets used for all numerical examples and the sizes of the datasets used.

\begin{table}
\centering
\footnotesize
\tabcolsep=0.175cm
\begin{tabular}{|c||r|r|r|r|r|r|r|r|r|}
\hline
\textbf{Spatial dimension}  &  \raisebox{\dimexpr-\height + 1.5ex\relax}{$\bpazocal{T}^1$}	& \raisebox{\dimexpr-\height + 1.5ex\relax}{$\bpazocal{T}^2$}  & \raisebox{\dimexpr-\height + 1.5ex\relax}{$\bpazocal{T}^3$}  & \raisebox{\dimexpr-\height + 1.5ex\relax}{$\bpazocal{T}^4$} & \raisebox{\dimexpr-\height + 1.5ex\relax}{$\bpazocal{T}^5$} & \raisebox{\dimexpr-\height + 1.5ex\relax}{$\bpazocal{T}^6$}  \\ \hline  \hline
\textbf{2} 					& $1,600$		& $6,241$		& $24,649$	& $97,969$	& $390,625$ 	& $1,560,001$    \\ \hline
\textbf{3} 					& $4,096$		& $29,791$	& $226,981$	& $1, 771,561$	& &   \\ \hline
\end{tabular}
\caption{Summary of the number of dofs associated with meshes of different resolutions.}
\label{table:mesh_hierarchies}
\end{table}

\begin{table}
\centering
\footnotesize
\begin{tabular}{|l||r|r|}
\hline
\multirow{2}{*}{\textbf{Example}}    &  \multicolumn{2}{c|}{\textbf{Branch network}}  \\ \cline{2-3} 
                		&  \multicolumn{1}{c|}{\textbf{Layers}} & \textbf{Act.}                               \\ \hline \hline
{Diff-TB}		& $2\times($Conv[40, 60, 100, 180] FFN[180, 80, 128]) & ReLU	 \\ \hline
{Diff-DP}		& $2\times($Conv[40, 60, 100, 180] FFN[180, 80, 128]) & ReLU	 \\ \hline
{Helm1D-TB} 	& FFN[1, 120, 120, 128]  & Tanh	 \\ \hline
{JumpDiff-TB} 	& FFN[1, 256, 256, 256, 128]  & Tanh	 \\ \hline
{Helm2D-TB}		& FFN[2, 256, 256, 256, 128]  & Tanh	 \\ \hline \hline
\multirow{2}{*}{\textbf{Example}}    &    \multicolumn{2}{c|}{\textbf{Trunk network}}    \\ \cline{2-3} 
                		&   \multicolumn{1}{c|}{  \textbf{Layers} } & \textbf{Act.}                                                                        \\ \hline 
{Diff-TB}		&  FFN[3, 80, 80, 128] &  Tanh \\ \hline
{Diff-DP}		&  FFN[3, 80, 80, 128] &  Tanh \\ \hline
{Helm1D-TB}	&  FFN[1, 150, 150, 128] &  Tanh \\ \hline
{JumpDiff-TB} 	&  FFN[2, 256, 256, 128] &  Tanh \\ \hline
{Helm2D-TB}		&  FFN[2, 256, 256, 128] &  Tanh \\ \hline
\end{tabular}
\caption{The summary of DeepONets' architectures.}
\label{tab:architecture}
\end{table}

\begin{table}[]
\centering
\footnotesize
\begin{tabular}{|l||r|r|r|r|r|r|r|r|r|}
\hline
\textbf{Example}    			& $\mathbf{N_S}$ 		& \textbf{Time  (mins)}   \\ \hline \hline
Diff-TB 				& 2,500		& $135.8$			   \\ \hline
Diff-DP 				& 2,500		& $136.2$			   \\ \hline
Diff-DP 				& 25,000		& $1,150.1$			   \\ \hline
Diff-DP 				& 100,000		& $4,972.8$			   \\ \hline
Helm1D-TB ($k_{\text{H}}=0$) 	 		& 2,500		& $9.8$			   \\ \hline
Helm1D-TB ($k_{\text{H}}=30$) 	 	& 2,500		& $10.1$			   \\ \hline
Helm1D-TB ($k_{\text{H}}=60$) 	 	& 2,500		& $12.5$			   \\ \hline
Helm1D-TB ($k_{\text{H}}=90$) 	 	& 2,500		& $18.3$			   \\ \hline
JumpDiff-TB 			& 2,500		& $121.4$			   \\ \hline
Helm2D-TB ($k_{\text{H}}=15$) 	 	& 2,500		& $156.2$			   \\ \hline
Helm2D-TB ($k_{\text{H}}=30$)  	& 2,500		& $892.4$			   \\ \hline
Helm2D-TB ($k_{\text{H}}=60$) 		& 2,500		& $5,123.8$		   \\ \hline
\end{tabular}
\caption{Summary of training information for all benchmark problems.}
\label{table:training_params}
\end{table}

\section{Trunk basis functions visualization}
\label{sec:basis_jump}
\begin{figure}[t]
\includegraphics[scale=0.5]{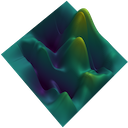}
\includegraphics[scale=0.5]{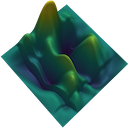}
\includegraphics[scale=0.5]{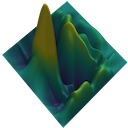}
\includegraphics[scale=0.5]{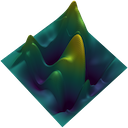}
\includegraphics[scale=0.5]{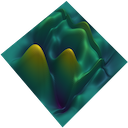}
\includegraphics[scale=0.5]{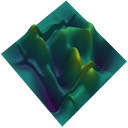}
\includegraphics[scale=0.5]{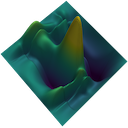}
\includegraphics[scale=0.5]{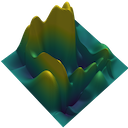}
\includegraphics[scale=0.5]{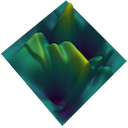}
\includegraphics[scale=0.5]{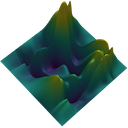}
\includegraphics[scale=0.5]{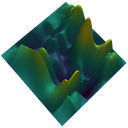}
\hspace{0.27cm}
\includegraphics[scale=0.5]{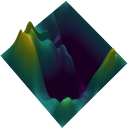}
\hspace{0.27cm}
\includegraphics[scale=0.5]{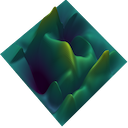}
\hspace{0.27cm}
\includegraphics[scale=0.5]{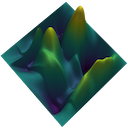}
\hspace{0.27cm}
\includegraphics[scale=0.5]{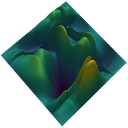}
\caption{An illustration of randomly selected $15$ trunk basis functions from the DeepONet ($p=128$) trained for \nameref{sec:jump_diff} test problem.}
\label{fig:tb_jump}
\end{figure}

\begin{figure}[t]
\includegraphics[scale=0.5]{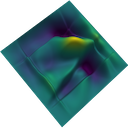}
\includegraphics[scale=0.5]{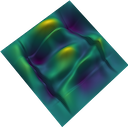}
\includegraphics[scale=0.5]{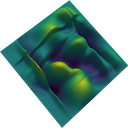}
\includegraphics[scale=0.5]{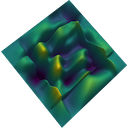}
\includegraphics[scale=0.5]{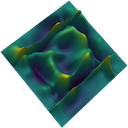}
\includegraphics[scale=0.5]{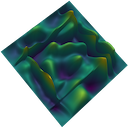}
\includegraphics[scale=0.5]{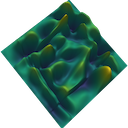}
\includegraphics[scale=0.5]{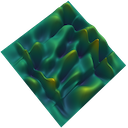}
\includegraphics[scale=0.5]{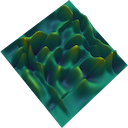}
\includegraphics[scale=0.5]{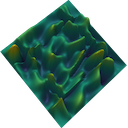}
\includegraphics[scale=0.5]{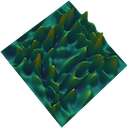}
\hspace{0.29cm}
\includegraphics[scale=0.5]{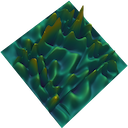}
\hspace{0.29cm}
\includegraphics[scale=0.5]{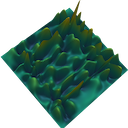}
\hspace{0.29cm}
\includegraphics[scale=0.49]{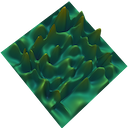}
\hspace{0.29cm}
\includegraphics[scale=0.49]{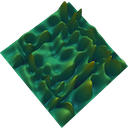}
\caption{An illustration of randomly selected $15$ trunk basis functions from the DeepONet ($p=128$) trained for \nameref{sec:jump_diff} test problem.
The basis functions are illustrated after performing QR decomposition.}
\label{fig:tb_jump_qr}
\end{figure}

\Cref{fig:tb_jump} and \Cref{fig:tb_jump_qr} visualize TB functions for \nameref{sec:jump_diff} benchmark problem before and after performing QR decomposition, respectively.

\section{Command line options}
\label{sec:command_line_options}
An example of command line options used to compose a DeepONet-based preconditioner can be found in \cref{fig:command_line_options}. 

\begin{figure}
\centering
\begin{tabular}{c}
  \begin{lstlisting}
-ksp_type gmres 
-pc_type mg 
-mg_levels_2_ksp_type richardson                            
-mg_levels_2_ksp_richardson_scale gamma_2
-mg_levels_2_pc_type jacobi
-mg_levels_1_ksp_type richardson                            
-mg_levels_1_pc_type composite
-mg_levels_1_pc_composite_type multiplicative
-mg_levels_1_pc_composite_pcs ksp,python,ksp
-mg_levels_1_sub_0_ksp_ksp_type richardson
-mg_levels_1_sub_0_ksp_ksp_max_it 1
-mg_levels_1_sub_0_ksp_pc_type jacobi
-mg_levels_1_sub_0_ksp_ksp_richardson_scale gamma_1
-mg_levels_1_sub_1_pc_python_type DONs.TB
-mg_levels_1_sub_1_pc_coarse_space_size  5
-mg_levels_1_sub_2_ksp_ksp_type richardson
-mg_levels_1_sub_2_ksp_ksp_max_it 1
-mg_levels_1_sub_2_ksp_pc_type jacobi
-mg_levels_1_sub_2_ksp_ksp_richardson_scale gamma_1
  \end{lstlisting}
  \end{tabular}
\caption{ An example of command line options used for composing DeepONet-based preconditioner in PETSc.
We precondition GMRES using three-level geometric MG with Jacobi smoother, which is multiplicatively composed with a DeepOnet {TB ($k=5$)} coarse-space (DONs.TB).
Variables gamma\_1 and gamma\_2 describe scaling parameters $\gamma_1, \gamma_2$ associated with the first and second levels, respectively.}
\label{fig:command_line_options}
\end{figure}

\section{Notation}
\label{sec:notations}
\Cref{tab:notation} summarizes notation used in the manuscript.

{
\footnotesize
\begin{longtable}{|p{1.3cm}|p{10.8cm}|}
\hline
\textbf{Symbol} 			&  \textbf{Description} 	\\ \hline \hline \endhead
$\Omega$  			& Computational domain \\ \hline
$\Gamma$  			& Dirichlet boundary \\ \hline
$d$  					& Spatial dimension of the problem \\ \hline
$\boldsymbol{\theta}$ 	& Parameters of the problem \\ \hline 
$\boldsymbol{\Theta}$ 	& Space of problem's parameters \\ \hline
$P$  					& Number of problem's parameters \\ \hline
$\pazocal{V}$  			& Hilbert space \\ \hline
\raisebox{\dimexpr-\height + 1.5ex\relax}{$\pazocal{V}'$}  			& Dual of Hilbert space \\ \hline
$\pazocal{V}_h$  		& The finite element space \\ \hline
$\pazocal{Y}$  			& Infinite-dimensional Banach space \\ \hline
$\pazocal{Y}_h$  		& Finite-dimensional Banach space \\ \hline
$u$  					& The solution of PDE \\ \hline
$\pazocal{A}$  			& Differential operator \\ \hline
$f$  					& Linear continuous form \\ \hline
\raisebox{\dimexpr-\height + 1.5ex\relax}{$a(\cdot, \cdot)$}  	& Parametrized bilinear form \\ \hline
\raisebox{\dimexpr-\height + 1.5ex\relax}{$f(\cdot)$}  		& Parametrized linear form \\ \hline
$\pazocal{T}$  			& Mesh \\ \hline
$L$  					& Number of meshes in the multilevel hierarchy \\ \hline
$h$	    				& Mesh size 				\\  \hline
$\psi$  				& Basis function of the FE space \\ \hline
$\Am$  				& Stiffness matrix 			\\  \hline
$\fv$  				& Right hand side 			\\  \hline
$\uv$  				& Solution 			\\  \hline
$n$	    				& Number of degrees of freedom (dofs) \\  \hline
$\pazocal{G}$	 		& Nonlinear mapping approximated by the DeepONet \\  \hline
$B$	    				&Branch network \\  \hline
$T$	    				&Trunk network \\  \hline
$p$	    				&Number of basis functions of the DeepONet \\  \hline
$nb_j$	    			& Number of points (sensor locations) used for discretization of $j$-th input function \\  \hline
$nf$	    				& Number of input function used in DeepONet \\  \hline
\raisebox{\dimexpr-\height + 1.5ex\relax}{$y^j$}	  	& Continuous representation of $j$-th branch input function \\  \hline
\raisebox{\dimexpr-\height + 1.5ex\relax}{$\qv^j$}	& Coordinate point used for evaluation of $j$-th input function \\ \hline
\raisebox{\dimexpr-\height + 1.5ex\relax}{$\yv^j$}	& Discrete representation of $j$-th branch input function \\  \hline
$\boldsymbol{\xi}$		& Coordinate point used for DeepONet inference \\  \hline
$\pazocal{D}$ 			& Dataset \\ \hline
$n_{\text{don}}$ 		& Number of points used for evaluation of solution for training the DeepONet \\ \hline
$N_s$ 				& Number of samples used for training the DeepONet \\ \hline
$\wv$ 				& Trainable parameters of DeepONet \\ \hline
$\pazocal{K}$ 			& Krylov subspace \\ \hline
$\pazocal{L}$ 			& Subspace orthogonal to residual \\ \hline
$\rv$ 				& Residual \\ \hline
$\Zm$ 				& Orthonormal matrix generated by the generalized Arnoldi algorithm \\ \hline
$\Vm$ 				& Basis of the Krylov space \\ \hline
$\Im$ 				& Identity matrix \\ \hline
\raisebox{\dimexpr-\height + 1.5ex\relax}{$\widebar{\Hm}$} & Upper Hessenberg matrix \\ \hline
$\ev$ 				& Canonical basis vector \\ \hline
$\zeta$ 				& Norm of the residual evaluated at initial iteration  \\ \hline
$m$ 					& Iteration number at which the restart procedure is invoked  \\ \hline
$\pv$ 				& Search direction \\ \hline
$\alpha$ 				& Step size used in CG algorithm  \\ \hline
$\beta$ 				& Conjugacy parameter used in CG algorithm  \\ \hline
$\Mm$ 				& Preconditioner 			\\  \hline
$\boldsymbol{\vartheta}$ 	& Iterate transformed by the preconditioner \\ \hline
$\Em$ 				& Error propagation operator \\ \hline
$\rho$ 				& Spectral radius \\ \hline
$\kappa$ 				& Condition number \\ \hline
$\Rm$ 				& Restriction operator \\ \hline
$\Pm$ 				& Prolongation operator \\ \hline
$\boldsymbol{\Pi}$ 		& Projection operator \\ \hline
$\Cm$ 				& Coarse space operator \\ \hline
$\Qm$ 				& Operator associated with an application of the coarse space step \\ \hline
$k$ 					& Dimension of coarse space \\ \hline
$S$ 					& Number of preconditioners used for construction of the composite preconditioner \\ \hline
$\gamma_s$ 			& Step size of $s$-th preconditioner \\ \hline
$\gamma_{k_{\text{H}}}$ 	& Optimal step size for Jacobi smoother for Helmholtz problem  \\ \hline
$\vv$ 				& Vector to which the preconditioner is applied \\ \hline
$b$ 					& Function used to impose Dirichlet boundary conditions in hard way \\ \hline
\raisebox{\dimexpr-\height + 1.5ex\relax}{$\widebar{\Qm}, \widebar{\Rm}$} & QR factors of tentative restriction operator \\ \hline
$\epsilon$ 			& Threshold value \\ \hline
$\Omega_s$  			& (Possibly) Overlapping $s$-th subdomain number \\ \hline
\raisebox{\dimexpr-\height + 1.5ex\relax}{$\widetilde{\Omega}_s$}  & Non-overlapping portion of the $s$-th subdomain \\ \hline
$\pazocal{I}$ & Index set of all dofs \\ \hline
\raisebox{\dimexpr-\height + 1.5ex\relax}{$\pazocal{I}_s$} & Index set of all dofs associated with the $s$-th subdomain \\ \hline
\raisebox{\dimexpr-\height + 1.5ex\relax}{$\widetilde{\pazocal{I}}_s$} & Index set of all dofs associated with the non-overlapping portion of the $s$-th subdomain \\ \hline
$n_s$ & Number of dofs associated with $s$-th subdomain \\  \hline
$K$ & Diffusion coefficient \\  \hline
$\ell, \sigma$ & Parameters specifying Gaussian random field \\  \hline
$\pazocal{U}$ & Uniform distribution \\  \hline
$\pazocal{N}$ & Gaussian normal distribution \\  \hline
$k_{\text{H}}, \lambda$ & Wave number and wavelength \\  \hline
$\Mm_{\mm}$ & Mass matrix \\  \hline
\caption{Summary of mathematical notations utilized in manuscript.}
\label{tab:notation}
\end{longtable}
}

\end{sloppypar}

{
\footnotesize
\section*{Acknowledgments}
{\footnotesize
We would like to thank Ansys, Inc.~for sponsoring this research.
AK gratefully acknowledges support of the Swiss National Science Foundation through the projects ``Multilevel training of DeepONets -- multiscale and multiphysics applications'' (206745), and  ``ML$^2$ -- Multilevel and Domain Decomposition Methods for Machine Learning'' (197041) as well as the  support of the Platform for Advanced Scientific Computing (PASC) under the project EXATRAIN.
GEK would like to acknowledge support by the DOE SEA-CROGS project (DE-SC0023191), the MURI-AFOSR FA9550-20-1-0358 project, and the ONR Vannevar Bush Faculty Fellowship (N00014-22-1-2795).}
{\footnotesize
AK would also like to extend her gratitude to Adar Kahana and Enrui Zhang for introducing her to HINTS and to Somdatta Goswami for generously sharing her extensive knowledge about DeepONets.
Moreover, she is thankful to Hardik Kothari for proofreading the manuscript and to Stefano Zampini for providing her with an introduction to petsc4py.}
}

{
\tiny
\bibliographystyle{siamplain}
\bibliography{biblio}
}

\end{document}